\documentclass[11pt,reqno]{amsart}
\usepackage{amsmath,mathtools}
\usepackage{graphicx}
\usepackage{epsfig,epsf,psfrag}
\usepackage{amsfonts}
\usepackage{setspace}
\usepackage{color}
\usepackage[font=small]{caption}
\usepackage[font=footnotesize]{subcaption}
\usepackage[top=1in, bottom=1in, left=1.2in, right=1.2in]{geometry}
\usepackage{afterpage}
\usepackage{bm}
\usepackage{multicol}
\usepackage{multirow} 

\usepackage{isomath}

\newcommand{\kp}{k_{\mathfrak p}}
\newcommand{\ks}{k_{\mathfrak s}}

\newcommand{\fp}{\mathfrak{p}}
\newcommand{\fs}{\mathfrak{s}}

\newcommand{\uinc}{\bm{u}^{inc}}
\newcommand{\usc}{\bm{u}^{sc}}

\newcommand{\bx}{\bm{x}}
\newcommand{\by}{\bm{y}}

\renewcommand{\phi}{\varphi}

\newtheorem{theorem}{Theorem}[section]
\newtheorem{corollary}[theorem]{Corollary}
\newtheorem{lemma}[theorem]{Lemma}

\newtheorem{remark}[theorem]{Remark}

\title[]{A robust and high precision algorithm for elastic scattering problems from cornered domains} 

\author{Jianan Yao}
\address{School of Mathematical Sciences, Zhejiang University,
	Hangzhou, Zhejiang 310027, China}
\email{yaojianan@zju.edu.cn}

\author{Jun Lai}
\address{School of Mathematical Sciences, Zhejiang University,
	Hangzhou, Zhejiang 310027, China}
\email{laijun6@zju.edu.cn}

\thanks{The work of JL was partially supported by the Key Project of Joint Funds for Regional Innovation and Development (No. U21A20425) and the “Xiaomi Young Scholars” program from Xiaomi Foundation.}



\date{\today}

\begin{document}

\maketitle 
\begin{abstract}
   The Navier equation is the governing equation of elastic waves, and computing its solution accurately and rapidly has a wide range of applications in geophysical exploration, material sciences, etc. In this paper, we focus on the efficient and high-precision numerical algorithm for the time harmonic elastic wave scattering problems from cornered domains via the boundary integral equations in two dimensions. The approach is based on the combination of Nystr\"om discretization, analytical singular integrals and kernel-splitting method, which results in a high order solver for smooth boundaries. It is then combined with the recursively compressed inverse preconditioning (RCIP) method to solve elastic scattering problems from cornered domains. Numerical experiments show the proposed approach is highly accurate with a stabilized error close to machine precision in various geometries. The algorithm is further applied to investigate the asymptotic behavior of density functions of boundary integral operators near corners. The results are highly consistent with the theoretical formula.
\end{abstract}


\section{Introduction}\label{introduction}

The analysis and computation of elastic scattering have received ever increasing
attention due to the  applications in many important areas such
as geological exploration, nondestructive testing, and material sciences \cite{ABG-15,Cakoni02,sneddon1980vd_uniqueness_3,lai2022fast_GGqApp2,laiFastInverseElastic2022}.  An accurate and efficient numerical method
plays an important role in many of these applications. For instance, in the inverse elastic obstacle scattering, often a large number of forward solving is required for methods based on the gradient descent~\cite{2018Inverse}. This paper is concerned with the scattering of a time-harmonic incident wave in a cornered domain filled with a homogeneous and isotropic elastic medium in two dimensions. In particular, we develop a robust and highly accurate numerical method for solving the elastic scattering problems from cornered domains under Dirichlet and Neumann boundary conditions. 

Cornered domains, characterized by sharp corners or edges in their boundaries, present unique challenges and opportunities for understanding the behavior of  elliptic and wave equations~\cite{Blsten2012CornersAS, GPcorner}. It is well known that the interplay between the geometry of the domain and the inherent properties of wave equations gives rise to fascinating scattering phenomena, such as corner singularities and diffraction effects~\cite{CharMi2016}. These singularities significantly impact the solutions of the governing equations and increase the difficulties of designing effective numerical schemes. Hence, although there are many computational methods available for the scattering problems of elastic waves~\cite{BXY2017, BU2014856, LELOUER20141, LY2022}, most of them focused on smooth domains. Among them, the method of boundary integral equations(BIE)~\cite{colton2013integral_potentialtheory} offers an attractive approach for solving the elastic boundary value problems. The advantages of BIE are that the discretization is only needed on the boundary of the domain when the medium is homogeneous, which reduces the dimension of computational domain by one, and the radiation condition at infinity is satisfied exactly for the exterior scattering problems. However, it is highly nontrivial to solve the BIE  to high precision in cornered domains due to the singular integrals and corner singularities. One usually needs to put a very dense mesh near the corner in order to fully resolve the singularities, which leads to a large, dense and ill-conditioned linear system. Direct inversion suffers from severe numerical loss. Over the years, several approaches have been proposed to overcome this difficulty, including the $L^2$ weighting method~\cite{bremer2012nystrom_corners}, corner rounding method~\cite{CharMi2016}, analytic expansion method~\cite{DensitySingular}, recursively compressed inverse preconditioning (RCIP)~\cite{helsing2013solving_RCIPtutorial},  etc. Among them, of particular interest is the method of RCIP, which was initially proposed in~\cite{helsing2008corner}, and was shown to be highly stable and accurate for solving the second kind boundary integral equations in cornered domains. However, they mostly focus on the Laplace and Helmholtz equations. Extension to elastic scattering is still challenging due to the complicated interconnection between compressional and shear waves~\cite{LL2019}. In particular, since the elastic Green's function is given by a second order tensor, it makes the decomposition of singular kernels in the elastic BIE highly nontrivial~\cite{Dong2020AHA}. Existing decomposition is often ineffective to be used in the high order discretization, especially for cornered domains. 

In this paper, we propose an effective solver for the elastic scattering problems from cornered domains. The method is based on the combination of Nystr\"om discretization, analytical singular integrals and kernel splitting method, which leads to a high order solver for smooth geometries.  When the boundary is piecewise smooth, we apply the non-uniform discretization with RCIP to fully resolve the corner singularities. In particular, we give a novel decomposition of the singular kernels for elastic integral boundary operators, which is highly compact and easy to compute, and show a detailed numerical implementation for the discretization of singular integrals. In the end, we demonstrate the effectiveness of the proposed approach by various numerical examples, and verify the leading order singularities for the density functions of BIE near the corner, which turns out to be highly consistent with the theoretical prediction.  

This paper is organized as follows. In Section \ref{problem_form}, we formulate the time harmonic elastic scattering problems in two dimensions and show the singular behavior of the solution for the elastic scattering near the corner. In Section \ref{potential_theory}, we introduce the elastic potential theory and derive the decompositions of kernel functions in elastic boundary integral operators based on the asymptotic analysis. In Section \ref{numerical_discretization}, a high order numerical discretization method for the elastic BIE in cornered domains is proposed, which combines the Nystr\"om discretization, singular integral quadrature and the recursively compressed inverse preconditioning. In Section \ref{numer_exp}, several numerical examples are given to illustrate the effectiveness of the proposed method. The paper is concluded in Section \ref{conclusion}.


\section{Problem formulation}\label{problem_form}

Consider a region $\Omega\subset \mathbb{R}^2$ with a piecewise smooth boundary $\Gamma$. We assume the number of corner points on $\Gamma$ is finite and $\Omega$ is locally diffeomorphic in the neighborhood of every corner point to an infinite cone with an opening angle $\theta_c\notin\{0,2\pi\}$. Denote $\bm{n}=(n_1, n_2)^\top$ and $\bm{\tau}=(\tau_1, \tau_2)^\top$ the unit
 normal and tangential vectors on $\Gamma$, respectively, where $\tau_1=-n_2,
 \tau_2=n_1$.  
 Assume $\Omega$ is filled with a homogeneous and isotropic elastic medium and illuminated by an elastic incident wave $\bm{u}^{inc}$. The total field $\bm{u}=(u_1,u_2)^\top$ satisfies the time-harmonic (with time dependence $e^{-i\omega t}$) Navier equation
\begin{equation}
	\mu \Delta \bm{u}+(\lambda+\mu) \nabla(\nabla \cdot \bm{u})+\rho \omega^{2} \bm{u}=0, \quad \bx\in \Omega, \label{navier}
\end{equation}
where $\mu$ and $\lambda$ are Lam\'e coefficients (with $\mu>0$ and $\lambda+\mu>0$), $\rho$ is the mass density, and $\omega$ is the angular frequency. Throughout, we assume all the material parameters are constant.

The total field $\bm{u}$ in equation \eqref{navier} is consisted of the incident field $\bm{u}^{inc}$ and the scattered field $\bm{u}^{sc}$.  There are typically two types of incident wave. One is the plane wave, which is a linear combination of the time-harmonic compressional plane wave $d\mathrm{e}^{\mathrm{i} \kp d\cdot	x}$ and the shear plane wave $d^{\perp}\mathrm{e}^{\mathrm{i}\ks d\cdot x}$, where $d=(\cos\theta, \sin\theta)^\top$ is
the unit propagation direction vector, $\theta\in [0, 2\pi)$ is the incident
angle, $d^{\perp}=(-\sin\theta, \cos\theta)^\top$ is an orthonormal vector of
$d$,
and
\[
k_{\fp}=\omega \sqrt{\frac{\rho}{\lambda+2 \mu}}, \quad k_{\fs}=\omega \sqrt{\frac{\rho}{\mu}},
\]
are the compressional  and shear wavenumber, respectively.  Another type of incident wave is the point source, which is given by equation \eqref{nav_green} in the next section. Both of them satisfy the Navier equation \eqref{navier} globally except at the source point.

Our main concern is to find the scattered field $\usc$ effectively from the cornered domains. The property of  $\usc$ depends on the boundedness of domain $\Omega$. Here we refer the elastic scattering problem as the interior problem if $\Omega$ is bounded, and call it the exterior problem if $\mathbb{R}^2\backslash \overline{\Omega}$ is bounded. 

\subsection{Interior elastic scattering problems}

We are concerned with two types of elastic scattering problems depending on the boundary condition of $\Gamma$.

\textit{Interior Dirichlet problems}: When the exterior of $\Omega$ is rigid, no deformation displacement occurs on the boundary $\Gamma$. Therefore, the total field $\bm{u}$ satisfies
\begin{align*}
	\bm{u}=0, \quad \bx \in \Gamma, 
\end{align*}
through which the scattered field $\usc$ satisfies the Dirichlet boundary value problem
\begin{align}
	\begin{aligned} \label{Dirichlet_problem}
		\begin{cases}
		\mu \Delta \usc+(\lambda+\mu) \nabla(\nabla \cdot \usc)+\rho \omega^{2} \usc = 0  , &\qquad \mbox{in} \quad \Omega,  \\
		\usc = -\uinc  ,&\qquad \mbox{on} \quad\Gamma.
		\end{cases}
	\end{aligned}
\end{align}

\textit{Interior Neumann problems}: 
Given a vector function $\bm{v}=(v_1, v_2)^\top$ and
a scalar function $v$, we define the scalar and vector curl operators by
\[
{\rm curl}\boldsymbol v=\partial_{x_1}v_2-\partial_{x_2}v_1, \quad
{\bf curl}v=(\partial_{x_2}v, -\partial_{x_1}v)^\top,
\]
respectively. The traction operator $\mathcal{T}$ on $\Gamma$ is defined by
\begin{align} \label{Traction1.4}
	\begin{aligned}
		\mathcal{T}\bm{u} &=2 \mu \frac{\partial \bm{u}}{\partial \bm{n}}+\lambda(\nabla \cdot \bm{u}) \bm{n}-\mu \bm{\tau}{\rm curl} \bm{u}. 
	\end{aligned}
\end{align}
We call that $\Gamma$ is traction free if
 \[\mathcal{T}\bm{u} = 0,\quad \bx\in \Gamma.\] 
In this case, the scattered field $\usc$ satisfies the Neumann boundary value problem  
\begin{align}
	\begin{aligned} \label{Neumann_problem}
		\begin{cases}
		\mu \Delta \usc+(\lambda+\mu) \nabla(\nabla \cdot \usc)+\rho \omega^{2} \usc = 0  , &\qquad \mbox{in} \quad \Omega,  \\
		\mathcal{T}\usc = -\mathcal{T}\uinc ,&\qquad \mbox{on} \quad \Gamma.
	\end{cases}
	\end{aligned}
\end{align}

\subsection{Exterior boundary value problems}
For the elastic scattering in an unbounded domain, one can first decompose the scattered field $\bm{u}^{sc}$ as
\begin{equation}
	\usc = \bm{u}_{\fp} + \bm{u}_{\fs},
\end{equation}
with 
\begin{equation}\label{helmdec}
	\bm{u}_{\fp}=-\frac{1}{k_{\fp}^{2}} \nabla(\nabla \cdot \usc), \quad \bm{u}_{\fs}=\frac{1}{k_{\fs}^{2}} \mathbf{curl}{\rm curl} \usc,
\end{equation}   
where $\bm{u}_{p}$ is called the compressional wave, and $\bm{u}_{s}$ is the shear wave. Equation \eqref{helmdec} is the Helmholtz decomposition of the scattered field $\usc$~\cite{LL2019}, as in the two dimensions both $\bm{u}_{\fp}$ and $\bm{u}_{\fs}$ satisfy the Helmholtz equation with wavenumber $k_\fp$ and $k_\fs$ respectively. The \textit{Exterior Dirichlet problems} and \textit{Exterior Neumann problems} are given the same as \eqref{Dirichlet_problem} and \eqref{Neumann_problem}, but coupled with the Kupradze-Sommerfeld radiation condition at infinity
\begin{equation} \label{Kupradze-Sommerfeld}
	\lim _{r \rightarrow \infty} r^{1 / 2}\left(\frac{\partial \bm{u}_{\fp}}{\partial r}-i k_{\fp} \bm{u}_{\fp}\right)=0 , \quad \lim _{r \rightarrow \infty} r^{1 / 2}\left(\frac{\partial \bm{u}_{\fs}}{\partial r}-i k_{\fs} \bm{u}_{\fs}\right)=0 , \quad r=|\bx|.
\end{equation}
 
\subsection{Well-posedness}
When the boundary $\Gamma$ is smooth, the existence and uniqueness for the solutions of the elastic boundary value problems have been well established in~\cite{sneddon1980vd_uniqueness_3}.
In particular, for the interior boundary value problems, when the wavenumber $k_\fp>0$ and $k_\fs> 0$, there exists a unique solution in the space of $(C^2(\Omega)\cap C^1(\overline{\Omega}))^2$ to equations \eqref{Dirichlet_problem} and \eqref{Neumann_problem} except a countable set of $k_\fp$ and $k_\fs$. These special wavenumbers are called the interior Dirichlet and Neumann eigenvalues, respectively.  For the exterior problems, as long as the wavenumbers $k_\fp$ and $k_\fs$ are positive, there always exists a unique solution  that satisfies the equations \eqref{Dirichlet_problem} and \eqref{Neumann_problem} coupled with the radiation condition \eqref{Kupradze-Sommerfeld}. Throughout this paper, unless otherwise stated, we only consider the regular wavenumbers by explicitly excluding these eigenvalues, so that the interior and exterior problems always admit a unique solution.

In case the boundary is piecewise smooth, due to the corner singularity, the solution for \eqref{Dirichlet_problem} and \eqref{Neumann_problem} only uniquely exists in some weighted Sobolev spaces~\cite{Cakoni02}. In particular, the behavior of the solution $\usc$ for the boundary value problems \eqref{Dirichlet_problem} or \eqref{Neumann_problem} in a neighborhood of the corner point can be described with the help of solutions of a homogeneous boundary value problem for the Lam\'e operator~\cite{hsiao2008boundary} in the infinite sector with the opening angle $\theta_c$. It was shown in~\cite{SRS1989} that the set of all solutions of this homogeneous problem has a basis consisting of the so-called “power solutions”
\begin{eqnarray}\label{powersolution}
	u(r,\theta) \sim \sum_{j=1}^{N}\sum_{k=1}^{m_j}\sum_{l=0}^{\kappa_{j,k}} a_{j,k,l}r^{\alpha_j} \sum_{s=0}^{l}(\log r)^s\omega_j^{l-s,k}(\theta) + O(r), 
\end{eqnarray}
where $(r,\theta)$ with $0\le \theta\le \theta_c$ are polar coordinates with origin at the corner point and $\{\omega_j^{l,k},1\le k\le m_j, 0\le l\le \kappa_{j,k}\}$ is a canonical system of Jordan chains of some linear operator pencil $A$~\cite{NA2022} corresponding to the eigenvalue $\alpha_j$. Here $m_j$ is the geometric multiplicity of $\alpha_j$ and $\kappa_{j,k}+1$ is the length of the $k$th Jordan chain. The number $N$ is chosen such that $0<{\rm Re}{(\alpha_j)}<1$. 

Due to the singularity at the corner, standard numerical method based on uniform mesh is usually inefficient and low order. Mesh refinement near the corner is unavoidable in numerical computation (except methods based on the analytical expansion~\cite{DensitySingular}), which on the other hand causes numerical instability and generates a large linear system. To address these issues, we propose a robust and high order numerical method to solve both the Dirichlet problems \eqref{Dirichlet_problem} and Neumann problems \eqref{Neumann_problem} in cornered domains based on the boundary integral equations.  

\section{Potential theory and asymptotic analysis}\label{potential_theory}
The boundary integral equation for elastic scattering is based on the potential theory of Navier equations. Given $\bx:=(x_1,x_2)$ and $\by:=(y_1, y_2)$ with $\bx \neq \by$, the free space Green's function for equation $\eqref{navier}$ in two dimensions is given by
\begin{equation}
	\mathbf{G}(\bx,\by) = \frac{i}{4\mu} \left\{ 
	H_0^{(1)} (k_\fs|\bx-\by|) \mathbf{I} 
	+ \frac{1}{k^2_\fs} \nabla_{\bx}\nabla_{\bx}  ^\top \left( H_0^{(1)} (k_\fs|\bx-\by|) - 
	H_0^{(1)} (k_\fp|\bx-\by|) \right) 
	\right\}, \label{nav_green}
\end{equation}
where $H_{0}^{(1)}$ is the first kind Hankel function of order zero, and $\mathbf{I}$ is the $2 \times 2$ identity matrix. 

When $\omega = 0$, equation $\eqref{navier}$ becomes the Lam\'e system\cite{hsiao2008boundary} and the corresponding Green's function is
\begin{align}\label{lam_green}
	\mathbf{E}(\bx,\by)=\frac{\lambda+3 \mu}{4 \pi \mu(\lambda+2 \mu)}\left\{\gamma(\bx, \by) \mathbf{I}+\frac{\lambda+\mu}{\lambda+3 \mu} \frac{1}{r^{2}}(\bx-\by)(\bx-\by)^{\top}\right\},
\end{align}
 with
\begin{align*}
	\gamma(\bx, \by)= -\ln|\bx-\by|.
\end{align*}
Let $\bm{\phi}=[\phi_1,\phi_2]^\top$ be a vector function defined on $\Gamma$. The single layer potential operator $S$ is defined by
	\begin{equation}\label{s_layer}
		S[\bm{\phi}](\bx)=\int_{\Gamma} \mathbf{G}(\bx,\by) 
		\bm{\phi}(\by) \textup{d} s_{\by}  , \quad \bx \notin \Gamma.
	\end{equation}
When $\bx$ approaches the boundary $\Gamma$, it becomes the single layer boundary operator
\begin{align}
	 	\mathcal{S}[\bm{\phi}](\bx)=\int_{\Gamma} \mathbf{G}(\bx,\by) 
	\bm{\phi}(\by) \textup{d} s_{\by}  , \quad \bx \in \Gamma. \label{single_layer_bd} 
\end{align}
The double layer potential operator $D$ is defined by
	\begin{equation}\label{d_layer}
		D[\bm{\phi}](\bx)=\int_{\Gamma} \mathbf{D}(\bx,\by) 
		\bm{\phi}(\by) \textup{d} s_{\by} , \quad \bx \notin \Gamma,
	\end{equation}
where the kernel is given by
	\begin{align*} 
		\mathbf{D}(\bx,\by)=\left(\mathcal{T}_{\by} \mathbf{G}(\bx,\by)\right)^{\top} 
	\end{align*}
with $\mathcal{T}_{\by}$ the traction operator acting on the variable $\by$. The corresponding boundary operator of double layer potential $D$ is 
 	\begin{align}
 	\mathcal{D}[\bm{\phi}](\bx)=\int_{\Gamma} \mathbf{D}(\bx,\by) 
 	\bm{\phi}(\by) \textup{d} s_{\by} , \quad \bx \in \Gamma.
 	\label{double_layer_bd}
 \end{align}
 The adjoint of double layer boundary operator is defined by
\begin{align}\label{double_layer}
	\mathcal{D}^{\prime}[\bm{\phi}](\bx)=\int_{\Gamma} \bm{\Sigma}(\bx, \by) \bm{\phi}(\by) \textup{d} s_{\by} , \quad \bx \in \Gamma,
\end{align}
with
\begin{align*} 
	\bm{\Sigma}(\bx,\by)=\mathcal{T}_{\bx} \mathbf{G}(\bx,\by).
\end{align*}
Note that both the double layer boundary operator and its adjoint are defined in the sense of Cauchy principal value. Similar to the acoustic case~\cite{colton2013integral_potentialtheory}, the double layer potential $D[\bm{\phi}](x)$ and the traction of single layer potential  $S[\bm{\phi}](x)$ will produce a jump when $\bx$ approaching the boundary $\Gamma$~\cite{hsiao2008boundary}, i.e.
\begin{align}\label{jumpcond}
	\begin{split}
	\lim\limits_{\epsilon \to 0^+} D[\bm{\phi}](\bx \pm \epsilon \bm{n})
	& = \pm \frac{1}{2} \bm{\phi}(\bx) + \mathcal{D}[\bm{\phi}](\bx), \\
	\lim\limits_{\epsilon \to 0^+} \mathcal{T}_{\bx}  S[\bm{\phi}](\bx \pm \epsilon \bm{n})
	& = \mp \frac{1}{2} \bm{\phi}(\bx) + \mathcal{D}^{\prime}[\bm{\phi}](\bx),
	\end{split}
\end{align}
where the unit normal vector $\bm{n}$ is towards the unbounded side of $\Gamma$. Note that \eqref{jumpcond} is only defined in the sense of almost everywhere as the limit at the corner point depends on the opening angle $\theta_c$~\cite{ATK97}. Based on \eqref{jumpcond}, we can formulate the boundary integral equations for elastic scattering in both Dirichlet and Neumann cases.

\subsection{Dirichlet problems}
For rigid boundary $\Gamma$, the solution $\usc_D$ to equation \eqref{Dirichlet_problem} can be represented as the double layer potential
\begin{align}
	\usc_D(\bx) = D[\bm{\phi}](\bx) =\int_{\Gamma}\left( \mathcal{T}_{\by} \bm{G}(\bx,\by) \right)^\top 
	\bm{\phi}_D(\by) \textup{d} s_{\by} , \quad \bx \in \Omega,
\end{align}
which results in a second kind boundary integral equation for the unknown vector function $\bm{\phi}_D$
\begin{align}\label{DND}
	\left(\pm \frac{1}{2}\mathcal{I} + \mathcal{D}\right)\bm{\phi}_D =  -\bm{u}^{inc} , \quad \bx \in \Gamma.
\end{align}
Here we take `$+$' for the exterior problem and `$-$' for the interior problem.
\begin{remark}
	Under the assumption that the wavenumber $k_\fp$ and $k_\fs$ are not the eigenvalue of the interior Dirichlet problems or Neumann problems, there exists a unique solution for the boundary integral equation \eqref{DND} in both interior and exterior cases~\cite{sneddon1980vd_uniqueness_3}.
\end{remark}
\begin{remark}
	For the exterior problems, one can also follow a similar approach in \cite{colton2013integral_potentialtheory} to formulate a combined layer potential to represent the scattered field
	\begin{align}\label{combinepot}
		\usc_D(x) = (D-iS)[\bm{\phi}_D](\bx) =\int_{\Gamma} (\mathcal{T}_{\by} \bm{G}-{i}\bm{G})(\bx,\by) 
		\bm{\phi}_D(\by) \textup{d} s_{\by} , \quad \bx \in \Omega,
	\end{align}
	which leads to the boundary integral equation
	\begin{align}\label{combinedIE}
			\left( \frac{1}{2}\mathcal{I} + \mathcal{D}-{ i}\mathcal{S}\right)\bm{\phi}_D =  -\bm{u}^{inc} , \quad \bx \in \Gamma.
	\end{align}
	It is resonance free for any positive wavenumber $k_\fp$ and $k_\fs$.
\end{remark}
\subsection{Neumann problems}
For the Neumann problem \eqref{Neumann_problem}, the scattered field $\usc_N$ can be represented by the single layer potential
\begin{align}
	\usc_N(\bx) = S[\bm{\phi}_N](\bx) =\int_{\Gamma} \bm{G}(\bx,\by) 
	\bm{\phi}_N(\by) \textup{d} s_{\by} , \quad \bx \in \Omega,
\end{align}
where the unknown vector function $\bm{\phi}_N$ satisfies the second kind  boundary integral equation 
\begin{align}\label{SNN}
	\left(\mp \frac{1}{2}\mathcal{I} + \mathcal{D}^{\prime}\right)\bm{\phi}_N =  -\mathcal{T}\bm{u}^{inc} , \quad \bx \in \Gamma.
\end{align}
Here we take '$-$' for the exterior problem and '$+$' for the interior problem. 
\begin{remark}
	Similar to the Dirichlet problem, there exists a unique solution to the boundary integral equation $\eqref{SNN}$ when $\kp$ and $\ks$ are not the interior eigenvalues.
\end{remark}
\begin{remark}	
	One can also represent the scattered field $\usc_N$ for the exterior Neumann problem by the combined layer potential as $\eqref{combinepot}$, so that the resulted boundary integral equation is resonance free for all the exterior Neumann problems. However, it will include the hyper-singular boundary operator, which makes the equation difficult to discretize. Discussion on the hyper-singular boundary integral formulation can be found in~\cite{BXY2017}.
\end{remark}

\subsection{Kernel splitting of the boundary operators}\label{kernel_split}
In order to construct a high order numerical method for the boundary integral equations \eqref{DND} and \eqref{SNN}, we need to analyze the singular kernels in the single and double layer boundary operators. To simplify the presentation, let us denote  $\bm{r} = \bx-\by, r = |\bx-\by|, \bm{r}^\bot=(-(x_2,-y_2),(x_1-y_1))^\top$, and  $\otimes$ the outer product of two vectors. We also set 
\begin{eqnarray}\label{defofL}
	\mathbf{L} = 	\begin{pmatrix}
		0 & 1 \\
		-1 & 0
	\end{pmatrix}.
\end{eqnarray}
Based on the identity~\cite[$\S$10.6.6]{olver2010nist} 
\begin{align}\label{hanklederi}
	\frac{\textup{d}}{\textup{d}z} \left( \frac{H_n^{(1)}(z)}{z^n} \right)
	= - \frac{H_{n+1}^{(1)}(z)}{z^n},
\end{align}
it yields
\begin{eqnarray}\label{derihankel}
	\begin{split}
	\partial_{x_i} H_0^{(1)}(k r) &=-k H_1^{(1)}(k r)\frac{x_i-y_i}{r}, \\
		\partial_{x_j}\partial_{x_i} H_0^{(1)}(k r) &=k^2 H_2^{(1)}(k r)\frac{(x_i-y_i)(x_j-y_j)}{r^2}-k H_1^{(1)}(k r)\frac{\delta_{ij}}{r}, \\
		\partial_{x_l}\partial_{x_j}\partial_{x_i} H_0^{(1)}(k r) &=-k^3 H_3^{(1)}(k r)\frac{(x_l-y_l)(x_i-y_i)(x_j-y_j)}{r^3}\\&\quad +k^2 H_2^{(1)}(k r)\frac{\delta_{ij}(x_l-y_l)+\delta_{il}(x_j-y_j)+\delta_{jl}(x_i-y_i)}{r^2},
	\end{split}
\end{eqnarray}
with $i,j,l=1$ or $2$. We thus get the explicit form for the kernel of single layer potential
\begin{align}
	\begin{aligned}
		\mathbf{G}(\bx,\by) 
		=&\frac{i}{4\mu}
		\left\{ H_0^{(1)}(k_\fs r) \mathbf{I} + \frac{1}{k_{\fs}^2}\left( k_{\fs}^{2} H_{2}^{(1)}\left(k_{\fs} r\right)-k_{\fp}^{2} H_{2}^{(1)}\left(k_{\fp} r\right) \right)  \frac{\bm{r}\otimes\bm{r}}{r^{2}} \right.\\
		&+ 
		\left.\frac{1}{k_\fs^2}\left( -k_{\fs} H_{1}^{(1)}\left(k_{\fs} r\right)+k_{\fp} H_{1}^{(1)}\left(k_{\fp} r\right) \right)  \frac{ \mathbf{I} }{r}
		\right\}.
	\end{aligned}
\end{align} 
 When $z\rightarrow 0$, it holds the following asymptotic properties for the Hankel functions ~\cite{olver2010nist}
\begin{align}
	\begin{split}\label{asyhankel}
		H_0^{(1)}(z) &= J_0(z) + i\left\{\frac{2}{\pi} \left( \ln\left(\frac{z}{2}\right) + \gamma \right) J_0(z) + O(z^2)\right\}, \\ 
		H_1^{(1)}(z) &= J_1(z) + i \left\{ \frac{2}{\pi} \left( -\frac{1}{z} + \ln\left(\frac{z}{2}\right)J_1(z) \right) - \frac{z}{2\pi}(1-2\gamma) + O(z^3)\right\}  ,
		\\
		H_2^{(1)}(z) &= J_2(z) + i \left\{ -\frac{4}{\pi}\frac{1}{z^2} - \frac{1}{\pi} + \frac{2}{\pi}\ln\left(\frac{z}{2}\right)J_2(z) - \frac{z^2}{4\pi}\left(\frac{3}{4}-\gamma\right) +O(z^4)\right\},
		\\
		H_3^{(1)}(z) &= J_3(z)+i \left\{ -\frac{8}{\pi} \left( \frac{2}{z^3} + \frac{1}{4z} + \frac{z}{32}\right) + \frac{2}{\pi}\ln\left(\frac{z}{2}\right)J_3(z) +O(z^3)\right\} ,
	\end{split}
\end{align}
where $\gamma$ is the Euler's constant and 
\begin{eqnarray}\label{asybessel}
 J_0(z)\sim 1+O(z^2),\quad J_1(z)\sim \frac{z}{2}+O(z^3), \quad J_2(z)\sim \frac{z^2}{8} +O(z^4), \quad J_3(z) \sim O(z^3).
\end{eqnarray} 
Note that all the remainders in equations \eqref{asyhankel} and \eqref{asybessel} are analytic with respect to $z\in \mathbb{C}$. We then obtain the asymptotic formula of $\mathbf{G}(\bx,\by)$ as $r\rightarrow 0$
\begin{align}\label{Gxy_leadingS}
	\mathbf{G}(\bx,\by) = \frac{1}{4\pi\mu} \left( -\frac{\lambda+3 \mu}{ (\lambda+2 \mu) }\ln(r) \mathbf{I} 
	+ \frac{\lambda + \mu}{\lambda + 2\mu} \frac{\bm{r} \otimes \bm{r}}{r^2}\right)+O(1)+O(r^2\ln r ), 
\end{align}
where the leading term is exactly $\mathbf{E}(x,y)$ for the Lam\'e equation. This explains that the single layer potential $S$ is continuous across the boundary and the corresponding boundary operator is compact. However, equation \eqref{Gxy_leadingS} is not convenient for high order discretization, we need another way to decompose $\mathbf{G}$.
\begin{lemma}
	It holds the decomposition that
	\begin{eqnarray}\label{decomp_G} 
		\mathbf{G}(\bx,\by) = -\frac{2}{\pi}\ln(r) {\rm Im}(\mathbf{G}(\bx,\by)) + \bm{R}_{\mathbf{G}}(\bx,\by),
	\end{eqnarray}
	where  ${\rm Im}(\mathbf{G}(\bx,\by))$ is the imaginary part of $\mathbf{G}(\bx,\by)$. Both ${\rm Im}(\mathbf{G}(\bx,\by))$ and the remainder
	$\bm{R}_{\mathbf{G}}(\bx,\by)$  are analytic on $\Gamma$ if $\Gamma$ is smooth. 
\end{lemma} 

It is obtained by the asymptotic expansions \eqref{asyhankel} and the identity that $H_n^{(1)}(z) = J_n(z)+iY_n(z).$ 
When $\bx\rightarrow\by$ on $\Gamma$, we can explicitly obtain
\begin{eqnarray}
	\begin{split}
		 \bm{R}_{\mathbf{G}}(\bx,\by)\rightarrow	&\frac{i}{4\mu}
			\left\{ \left(1 + \frac{2 i}{\pi} \left( \ln\left(\frac{\ks}{2}\right) + \gamma \right) \right) \mathbf{I} - \frac{1}{k_{\fs}^2} \frac{i}{\pi}(k_{\fs}^{2}  -  k_{\fp}^{2}) \mathbf{M}\right.\\
			&+ 
			\left.\frac{1}{k_\fs^2}\left( -k^2_{\fs}\left(\frac{1}{2} + i \left( \frac{1}{\pi}   \ln\left(\frac{\ks}{2}\right)  - \frac{1}{2\pi}(1-2\gamma) \right) \right)\right.\right.\\&\left.\left.+k^2_{\fp} \left(\frac{1}{2} + i \left( \frac{1}{\pi}   \ln\left(\frac{\kp}{2}\right)  - \frac{1}{2\pi}(1-2\gamma) \right) \right) \right) \mathbf{I} 
			\right\}.
		\end{split}
\end{eqnarray}
where $\mathbf{M}=\lim_{\bx\rightarrow\by}\frac{\bm{r}\otimes\bm{r}}{r^{2}}$, which can be easily evaluated once the parametrization of $\Gamma$ is given. 

Similarly, by making use of equation \eqref{hanklederi}, the kernel for the double layer boundary operator can be decomposed as follows.
\begin{theorem}\label{decomp_D_thm}
	It holds the decomposition that
	\begin{align}\label{decomp_D}
		\mathbf{D}(\bx,\by)&=\left(\mathcal{T}_{\by} \mathbf{G}(\bx,\by)\right)^{\top}	=\left(\mathcal{T}_{\by} \mathbf{E}(\bx,\by)\right)^{\top} 
		+ \ln(r)\mathbf{F}_{\by}(\bx,\by) + \bm{R}_{\mathbf{D}}(\bx,\by).
	\end{align}
	where 
	\begin{align}
		\mathcal{T}_{\by} \mathbf{E}(\bx, \by) &=\frac{\bm{n}(\by) \cdot \bm{r}}{2 \pi r^{2}} \mathbf{I}-\frac{\mu}{\lambda+2 \mu} \frac{\bm{\tau}(\by) \cdot \bm{r}}{2 \pi r^{2}}
	\mathbf{L}  \notag \\
		&\quad -\frac{\lambda+\mu}{\lambda+2 \mu}\left(-\frac{(\bm{n}(\by) \cdot \bm{r}) \bm{r} \otimes \bm{r}}{\pi r^{4}}+\frac{(\bm{n}(\by) \cdot \bm{r})}{2 \pi r^{2}} \mathbf{I}\right),	\label{TyE}\\
		\mathbf{F}_{\by}(\bx,\by) &= -\frac{2}{\pi} \ \textup{Im} \left\{ \left(\mathcal{T}_{\by} \mathbf{G}(\bx,\by)\right)^{\top} \right\}, 
	\end{align}
	and $\bm{R}_{\mathbf{D}}(\bx,\by)\rightarrow 0$ when $\bx\rightarrow\by$.
	If the boundary $\Gamma$ is smooth, both $\mathbf{F}_{\by}(\bx,\by)$ and the remainder term $\bm{R}_{\mathbf{D}}(\bx,\by)$ in $\eqref{decomp_D}$ are analytic on $\Gamma$.   
\end{theorem}
\begin{remark}
	Here $\left(\mathcal{T}_{\by} \mathbf{E}(\bx,\by)\right)^{\top}$ is exactly the kernel of double layer potential of Lam\'e equation. The second term in $\left(\mathcal{T}_{\by} \mathbf{E}(\bx,\by)\right)^{\top}$ is of order $O\left(r^{-1}\right)$ as $r=|x-y| \rightarrow 0,$, which makes the boundary integral existed in the sense of Cauchy principal value~\cite{hsiao2008boundary}.	
\end{remark}

\begin{proof}
By equation \eqref{hanklederi} and the definition of traction operator $\mathcal{T}$ in \eqref{Traction1.4}, we have
\begin{align}
	\begin{aligned}
\mathcal{T}_{\by} \mathbf{G}(\bx,\by)
		&=\frac{i}{4}
		\left\{\frac{k_{\fs} H_{1}^{(1)}\left(k_{\fs} r\right)}{r} \bm{L}_1(\bx,\by)		
		+\frac{\lambda}{\lambda+2 \mu} \frac{k_{\fp} 
			H_{1}^{(1)}\left(k_{\fp} r\right)}{r}\bm{L}_2(\bx,\by)\right.\\
		&\quad+\frac{2}{k_{\fs}^{2}}\left(\frac{\left(k_{\fs}^{3} H_{3}^{(1)}\left(k_{\fs} r\right)-k_{\fp}^{3} H_{3}^{(1)}\left(k_{\fp} r\right)\right)}{r^{3}} \bm{L}_3(\bx,\by)\right.\\&\quad\left.
		\left.-\frac{\left(k_{\fs}^{2} H_{2}^{(1)}\left(k_{\fs} r\right)-k_{\fp}^{2} H_{2}^{(1)}\left(k_{\fp} r\right)\right)}{r^{2}} \bm{L}_4(\bx,\by)\right)\right\} ,
	\end{aligned}
\end{align}
where \begin{eqnarray}
	\begin{split}
		\bm{L}_1(\bx,\by) &= 2(\bm{n}(\by) \cdot \bm{r}) \mathbf{I}-\bm{\tau}(\by) \otimes \bm{r}^{\perp}, \quad
		\bm{L}_2(\bx,\by)=  \bm{n}(\by) \otimes \bm{r}, \\
		\bm{L}_3(\bx,\by)&= (\bm{n}(\by)\cdot \bm{r}) \bm{r}\otimes\bm{r},\quad 
		\bm{L}_4(\bx,\by)=(\bm{n}(\by) \cdot \bm{r}) \mathbf{I}+\bm{n}(\by) \otimes \bm{r}+\bm{r} \otimes \bm{n}(\by).
	\end{split}
\end{eqnarray}
Using the expansion \eqref{asyhankel}, the asymptotic expansion for the double layer boundary operator can be rewritten as
\begin{align}\label{TyGexpan}
	\begin{aligned}
\mathcal{T}_{\by} \mathbf{G}(\bx,\by)
		&=\frac{\bm{L}_1(\bx,\by)}{2 \pi r^2}
		- k_\fs\ln(r)J_1(k_\fs r)\frac{\bm{L}_1(\bx,\by)}{2\pi r} 
		\\&\quad
		+\frac{\lambda}{\lambda+2 \mu} \frac{\bm{L}_2(\bx,\by)}{2 \pi r^{2}} - \frac{\lambda}{\lambda+2 \mu}
		\frac{k_\fp\ln(r)J_1(k_\fp r)}{2\pi r}\bm{L}_2(\bx,\by)
		\\&\quad
		-\frac{2}{k_{\fs}^{2}} \left\{ -\left(k_{\fs}^{2}-k_{\fp}^{2}\right) \frac{\bm{L}_3(\bx,\by)}{2 \pi r^{4}} + \frac{k_\fs^4-k_\fp^4}{8} \frac{\bm{L}_3(\bx,\by)}{2 \pi r^{2}} \right\} 
		\\&\quad
		-\frac{2}{k_{\fs}^{2}} \left( k_\fs^3\ln(r)J_3(k_\fs r)-k_{\fp}^3\ln(r)J_3(k_\fp r) \right)
		\frac{\bm{L}_3(\bx,\by)}{2\pi r^{3}}
		\\&\quad
		-\frac{\left(k_{\fs}^{2}-k_{\fp}^{2}\right)}{k_{\fs}^{2}} \frac{\bm{L}_4(\bx,\by)}{2 \pi r^{2}}
		+\frac{2}{k_{\fs}^{2}} \left( k_\fs^2\ln(r)J_2(k_\fs r)-k_{\fp}^2\ln(r)J_2(k_\fp r) \right)
		\frac{\bm{L}_4(\bx,\by)}{2 \pi r^{2}} \\&\quad
		+O(\bm{r}).
	\end{aligned}
\end{align} 
Merging the logarithmic singular part and using the fact that
\begin{align}
	&\bm{r} = (\bm{n}(\by)\cdot\bm{r})\bm{n}(\by) + (\bm{\tau}(\by)\cdot\bm{r})\bm{\tau}(\by), 
	\quad
	\bm{r}^\bot = -(\bm{\tau}(\by)\cdot\bm{r})\bm{n}(\by) + (\bm{n}(\by)\cdot\bm{r})\bm{\tau}(\by),  
	\\
	&\textbf{I} = \bm{n}(\by)\otimes\bm{n}(\by) + \bm{\tau}(\by)\otimes\bm{\tau}(\by) , 
	\quad
	\mathbf{L} = \bm{n}(\by)\otimes\bm{\tau}(\by) - \bm{\tau}(\by)\otimes\bm{n}(\by),
\end{align}
where $\mathbf{L}$ is defined in \eqref{defofL}, we obtain the decomposition \eqref{decomp_D}.
\end{proof}

One can analogously derive the asymptotic expansion for the kernel of the adjoint of the double layer boundary operator $\mathcal{D}'$.
\begin{corollary} It holds the decomposition that
\begin{align}\label{decomp_S}
	\bm{\Sigma}(\bx,\by)=\mathcal{T}_{\bx} \mathbf{G}(\bx,\by)
	= \mathcal{T}_{\bx} \mathbf{E}(\bx, \by) + \ln(r)\mathbf{F}_{\bx}(\bx,\by) + \bm{R}_{\mathbf{\Sigma}}(\bx,\by),
\end{align}
with 
\begin{align}
	\mathcal{T}_{\bx} \mathbf{E}(\bx, \by) &=-\frac{\bm{n}(\bx) \cdot \bm{r}}{2 \pi r^{2}} \mathbf{I}+\frac{\mu}{\lambda+2 \mu} \frac{\bm{\tau}(\bx) \cdot \bm{r}}{2 \pi r^{2}}
\mathbf{L} \notag \\
	&+\frac{\lambda+\mu}{\lambda+2 \mu}\left(-\frac{(\bm{n}(\bx) \cdot \bm{r}) \bm{r} \times \bm{r}}{\pi r^{4}}+\frac{(\bm{n}(\bx) \cdot \bm{r})}{2 \pi r^{2}} \mathbf{I}\right),\label{TxE} \\
	\mathbf{F}_x(\bx,\by) &= -\frac{2}{\pi} \ \textup{Im} \left\{ \mathcal{T}_{\bx} \mathbf{G}(\bx,\by) \right\}, 
\end{align}
and $\bm{R}_{\mathbf{\Sigma}}(\bx,\by)\rightarrow 0$ when $\bx\rightarrow\by$. Both $\mathbf{F}_{\bx}(\bx,\by)$ and  $\bm{R}_{\mathbf{\Sigma}}(\bx,\by)$ are analytic on $\Gamma$ if $\Gamma$ is smooth.	
\end{corollary}

 We conclude from \eqref{Gxy_leadingS}, \eqref{decomp_D} and \eqref{decomp_S} that the principal parts of the kernel functions of elastic wave equation are the kernel functions of Lam\'e equation, which is similar to the relation between Helmholtz equation and Laplace equation. The difference is that the double layer boundary operator and its adjoint in elastic equations have Cauchy singularity even after the extraction of identity operator. Conventional quadrature method based on weakly singular integrals can not give a high order discretization. In the next section, we will make use of the kernel splitting results \eqref{decomp_G}, \eqref{decomp_D} and \eqref{decomp_S} to obtain a high order Nystr\"om discretization for the singular boundary operators.

\section{Numerical discretization}\label{numerical_discretization}

Nystr\"om  discretization is a commonly used numerical approach for boundary integral equations~\cite{ATK97,kress1989linear}. For smooth geometries, Nystr\"om method based on trigonometric interpolation seems to be the most effective way to discretize singular integrals and converges exponentially fast~\cite{Dong2020AHA}. On the other hand, it is less effective in the geometries with corners as extremely many Fourier modes are needed to fully resolve the corner singularity. Here we choose to discretize the boundary integral equation based on the piecewise polynomial interpolation. In particular, we make use of the composite Gauss-Legendre nodes and the analytical singular integrals to discretize the kernels given by equations \eqref{Gxy_leadingS}, \eqref{decomp_D}, and \eqref{decomp_S}. We will show that the singular quadrature based on such kernel-splitting method with uniform partition can achieve 15-digits or more accuracy on the smooth geometries. However, these digits can be easily lost in the cornered domains as uniform mesh is difficult to fully resolve the corner singularities. In order to construct a high order and efficient numerical method for the elastic scattering problems from cornered domains, we adopt a graded mesh near the corners to discretize the boundary. Through this we can resolve the corner singularity but also increase the size of the discretized matrix and the condition number. Directly inverting the resulted linear system is numerically unstable due to the round-off errors of small weights~\cite{bremer2012nystrom_corners}.  One therefore needs a method to keep the discretized linear system well-conditioned with a relatively low computational cost. For that purpose, the recursively compressed inverse preconditioning method (RCIP)~\cite{helsing2013solving_RCIPtutorial} is employed for the discretized scattering problems from cornered domains. 

%

In the following, we will give the details of the Nystr\"om method for the elastic scattering using the composite Gauss-Legendre quadrature and analytical singular integrals, and then briefly describe the RCIP method. 

\subsection{Discretization of singular integrals}
\begin{figure}[h]
	\centering
	\includegraphics[width=0.5\textwidth]{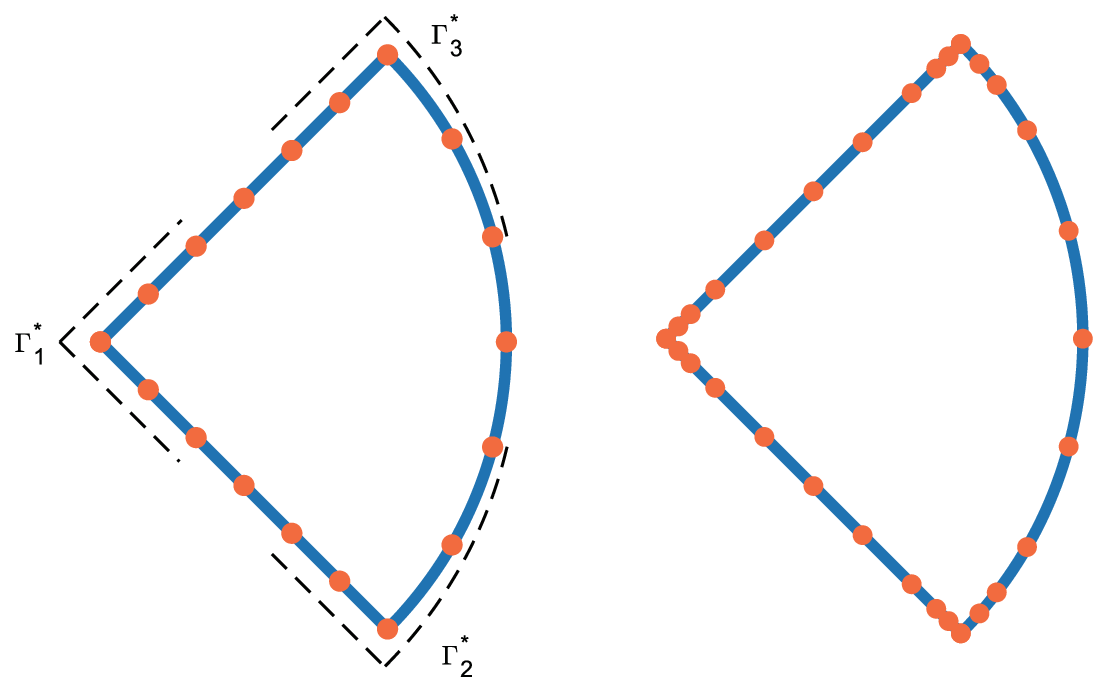}
	\vspace{-2pt}
	\caption{Discretization of a piecewise smooth boundary $\Gamma$. On the left is a coarse mesh, where each smooth component of $\Gamma$ is divided equally into $N$ panels.  For each corner point, we use the nearest four coarse panels to represent its neighborhood. On the right is the refined mesh obtained by dyadically subdividing the panels adjacent to the corners.}
	\label{meshes}
\end{figure}

 Assume that the piecewise smooth boundary $\Gamma$ is parameterized  by
\begin{align*}
	\bx(s) = (x_1(s),x_2(s)), \quad s \in[0,L].
\end{align*}
We first divide each smooth component of $\Gamma$ into $N$ uniform panels, where $N$ can be different for each component.  We call the uniform mesh as the coarse mesh. Due to the corner singularity, we apply the non-uniform discretization in the neighborhood of each corner point using graded mesh. More specifically, assume the first panel is $[0,1]$ with $0$ corresponding to the corner point, we subdivide it into $n_{sub}+1$ smaller panels as $[0,(1/2)^{n_{sub}}]\cup [(1/2)^{n_{sub}},(1/2)^{n_{sub}-1}]\cup \dots\cup [1/2,1]$. See Figure \ref{meshes} for an illustration of the coarse and refined meshes.

Denote $\Gamma_i$ the $i$th panel, $ i=1, 2,\cdots,N$. For the panel $\Gamma_i$, we further discretize it by  Gauss-Legendre quadrature with nodes $s_{ij}$ and weights $w_{ij}$, $j=0,1,\cdots,n-1$. Here  denote the corresponding points on the panel $\Gamma_i$ by $\bx_{ij} = \bx(s_{ij})$ or $\by_{ij} = \by(s_{ij})$.  Given $\bx \in \Gamma_l$, $l=1,2,\cdots,N$, consider the boundary integral in the form of
\begin{equation}\label{S_iintegral}
	\mathcal{K}_i[\bm{\phi}](\bx) = \int_{\Gamma_i} \bm{K}(\bx, \by) 
	\bm{\phi}(\by) \textup{d} s_{\by}, \quad \bx \in \Gamma_l,  
\end{equation}
where $\bm{\bm{K}}(\bx,\by)$ is the kernel function with singularity at $\bx=\by$. When $\Gamma_k$ and $\Gamma_i$ are distant from each other, the usual Nystr\"om discretization reads
\begin{equation}\label{Nystrom}
	\int_{\Gamma_i} \bm{K}(\bx, \by) 
	\bm{\phi}(\by) \textup{d} s_{\by} \approx \sum_{j=0}^{n-1} \bm{K}(\bx, \by_{ij}) 
	\bm{\phi}(\by_{ij})w_{ij}J_{i}(s_{ij})  
\end{equation}
where $J_{i}(s)$ is the Jacobian on panel $\Gamma_i$. It will achieve high order accuracy except for adjacent or self-interacting panels of $\Gamma_l$, i.e. $l=i-1, i, i+1$, due to the singularity of the kernel function. One remedy is to switch the quadrature to the generalized Gaussian quadrature scheme~\cite{bremer2010nonlinear_GGq}. It is accurate for logarithmic singularity but less effective for Cauchy singularity. In this paper, we make use of the decomposition \eqref{decomp_G}, \eqref{decomp_D}, and \eqref{decomp_S} to split each kernel into a combination of an analytically integrable singular part and a smooth non-singular part. For the singular part, one can integrate it by the analytical formula, while for the smooth part, the usual Gaussian quadrature \eqref{Nystrom} will give a high order accuracy. More specifically, we first transform the integral on $\Gamma_i$ to $[-1,1]$ by a suitable change of variable. Substituting the parametrization $\bx(t)$ and $\by(s)$ into the expressions \eqref{decomp_G}, \eqref{decomp_D}, and \eqref{decomp_S}, one can convert the kernel functions into the form of
\begin{eqnarray}\label{splitform}
	\bm{K}(\bx(t), \by(s)) =\ln(|t-s|)\bm{K}_1(t,s)+\frac{1}{t-s}\bm{K}_2(t,s)+\bm{K}_3(t,s),
\end{eqnarray}
 which swaps the singularity from $\bx$ and $\by$ to the parameters $t$ and $s$, and $\bm{K}_i,i=1,2,3$ are all analytical functions on $\Gamma_i$. For instance,  we can use \eqref{decomp_D} to rewrite the kernel $\bm{D}(\bx,\by)$ as
\begin{eqnarray}
	\begin{split}
	&\mathbf{D}(\bx(t),\by(s))	\\= &\ln(|t-s|)\mathbf{F}_{\by}(\bx,\by) +\frac{1}{t-s}\left(-\frac{(t-s)\mu}{\lambda+2 \mu} \frac{\bm{\tau}(\bx) \cdot \bm{r}}{2 \pi r^{2}}
	\mathbf{L}^{\top}\right) \\&+
	\left((\ln(r)-\ln(|t-s|))\mathbf{F}_{\by}(\bx,\by)+\left(\mathcal{T}_{\by} \mathbf{E}(\bx,\by)\right)^{\top}  +\frac{\mu}{\lambda+2 \mu} \frac{\bm{\tau}(\bx) \cdot \bm{r}}{2 \pi r^{2}}
	\mathbf{L}^{\top} + \bm{R}_{\mathbf{D}}(\bx,\by)\right).		
	\end{split}
\end{eqnarray}

\begin{remark}
	In theory, when $\bx\ne \by$, the value of $\bm{R}_{\mathbf{G}}$, $\bm{R}_{\mathbf{D}}$, and $\bm{R}_{\mathbf{\Sigma}}$ can be obtained by taking the difference of the left and right hand sides of equations \eqref{decomp_G}, \eqref{decomp_D} and \eqref{decomp_S}, respectively. However, in the numerical implementation, when $\bx$ is very close to $\by$, one should switch to the asymptotic expansions of $\bm{R}_{\mathbf{G}}$, $\bm{R}_{\mathbf{D}}$, and $\bm{R}_{\mathbf{\Sigma}}$ based on equations \eqref{asyhankel} and \eqref{asybessel} to avoid numerical round-off errors.	 
\end{remark}

Then the integral \eqref{S_iintegral} can be approximated by
\begin{align}\label{singular_int}
	\begin{split}
	&\int_{\Gamma_i} \bm{K}(\bx, \by) 
	\bm{\phi}(\by) \textup{d} s(\by)\\=&\int_{-1}^{1}\left[\ln(|t-s|)\bm{K}_1(t,s)+\frac{1}{t-s}\bm{K}_2(t,s)+\bm{K}_3(t,s)\right]\bm{\phi}(s)J_i(s)ds\\
	\approx &\int_{-1}^{1} \ln(|t-s|) \sum_{j=0}^{n-1} c_{ij} P_j(s)+\frac{1}{t-s}\sum_{j=0}^{n-1} d_{ij}P_j(s) + \sum_{j=0}^{n-1} f_{ij}P_j(s) ds\\
	=& \sum_{j=0}^{n-1} \left( c_{ij}\int_{-1}^{1} \ln(|t-s|) P_j(s)ds+d_{ij}\int_{-1}^{1}\frac{1}{t-s} P_j(s)ds + f_{ij}\int_{-1}^{1} P_j(s) ds\right)
	\end{split}
\end{align}
where $P_j(s)$ is the Legendre polynomial of degree $j$, $c_{ij}, d_{ij}$, and  $f_{ij}$ are the expansion coefficients of Legendre polynomials for $\bm{K}_1\bm{\phi}J_i, \bm{K}_2\bm{\phi}J_i$ and $\bm{K}_3\bm{\phi}J_i$. These coefficients can be obtained by using the $n\times n$ interpolation matrix $U=(u_{jk})$ that maps the value at the Guass-Legendre nodes to the coefficients of $P_j$. In other words, it holds
\begin{eqnarray}\label{cdf}
	c_{ij} = \sum_{k=0}^{n-1}u_{jk} [\bm{K}_1\bm{\phi}J_i](s_{ik}), \quad 
	d_{ij} = \sum_{k=0}^{n-1}u_{jk} [\bm{K}_2\bm{\phi}J_i](s_{ik}),\quad
	f_{ij} = \sum_{k=0}^{n-1}u_{jk} [\bm{K}_3\bm{\phi}J_i](s_{ik}).
\end{eqnarray}
Denote
\begin{eqnarray}\label{Asinteg}
	\begin{split}
	L_j &= \int_{-1}^{1} \ln(|t-s|) P_j(s)ds, \quad 
	C_j &= \int_{-1}^{1} \frac{1}{t-s} P_j(s)ds, \quad
	I_j &= \int_{-1}^{1}P_j(s)ds.
\end{split}
\end{eqnarray}
The value of $I_j$ can be easily obtained by the orthogonality of Legendre polynomials. The evaluation of $L_j$ and $C_j$ will be discussed in the next subsection. Once they are found, by combining \eqref{cdf} and \eqref{Asinteg}, we can rewrite the integral \eqref{singular_int} as
\begin{eqnarray}
	\begin{split}
	\int_{\Gamma_i} \bm{K}(\bx, \by) 
	\bm{\phi}(\by) \textup{d} s(\by)&\approx \sum_{j=0}^{n-1}\left( \sum_{k=0}^{n-1}u_{jk} \left(L_j\bm{K}_1\bm{\phi}J_i+ C_j\bm{K}_2\bm{\phi}J_i+  I_j\bm{K}_3\bm{\phi}J_i\right)
	\right)\\
	& = \sum_{k=0}^{n-1}\left[ \sum_{j=0}^{n-1} u_{jk} \left(L_j\bm{K}_1+ C_j\bm{K}_2+  I_j\bm{K}_3\right)\right]\bm{\phi}J_i(s_{ik}) \\
	& = \sum_{k=0}^{n-1}\tilde{w}_{ik}\bm{\phi}J_i(s_{ik}),
\end{split}
\end{eqnarray}
where $\tilde{w}_{ik} = \sum_{j=0}^{n-1} u_{jk} \left(L_j\bm{K}_1+ C_j\bm{K}_2+  I_j\bm{K}_3\right)J_i$ can be taken as the modified singular quadrature weights for the integral \eqref{S_iintegral}. Note that they not only depend on the Jacobian $J_i$ at $s_{ik}$, but also the target and source points $\bx(t)$ and $\by(s)$.

\subsection{Analytical singular integrals} \label{subCorrected singular integral}
Through the discussion above, the discretization of singular integrals is transformed into the evaluation of \eqref{Asinteg}. In this part, we focus on the computation of $L_j$ and $C_j$, $j=0,\cdots,n-1$ and note that they only need to be computed once. It is worth mentioning that the singular integrals for monomials have been discussed in~\cite{helsing2009integral_kpinitial}. We proceed by evaluating  $C_j$ first, as the evaluation of $L_j$ is based on $C_j$. 
\subsubsection{Cauchy singular integrals}
For a fixed $t \in \mathbb{R} $ and $t\ne \pm 1$, it holds
\begin{align*}
		C_0&=\int_{-1}^1 \frac{1}{t-s} P_0(s) \textup{d} s=-\ln \left|\frac{t-1}{t+1}\right|, \\
		C_1&=\int_{-1}^1 \frac{1}{t-s} P_1(s) \textup{d} s=-2+t C_0,
\end{align*}
where the integral exists in the sense of Cauchy principal value when  $t \in (-1,1)$. 
For $j \geq 2$, we can make use of the recursive formula for Legendre polynomials
\begin{equation*}
	jP_j(s) = (2j-1)sP_{j-1}(s)-(j-1)P_{j-2}(s),
\end{equation*}
to obtain the recursion formula of $C_j$,
\begin{align*}
	jC_j = (2j-1)tC_{j-1} - (j-1)C_{j-2}.
\end{align*}

\subsubsection{Logarithmic singular integrals}
We will make use of the Cauchy singular integral $C_j$ to compute the logarithmic integral $L_j$. More specifically, it holds
\begin{align*}
	L_0=&\int_{-1}^1 \ln |t-s| P_0(s) \textup{d} s=(t+1) \ln |t+1|-(t-1) \ln |t-1|-2, \\
	L_1=&\int_{-1}^1 \ln |t-s| P_1(s) \textup{d} s=t L_0-\left(\frac{(t+1)^2}{2} \ln |t+1|-\frac{(t-1)^2}{2} \ln |t-1|-t\right).
\end{align*}
For $j \geq 2$, according to the formula
\begin{equation*}
	\frac{\textup{d}\left(P_{j+1}(s)-P_{j-1}(s)\right)}{\textup{d} s} = (2 j+1) P_j(s),
\end{equation*}
we have
\begin{align*}
	\begin{aligned}
		j L_j =& j \int_{-1}^1 \ln |t-s| P_j(s) \textup{d} s \\
		=&\left\{\begin{array}{l}
			-\left[L_j-L_{j-2}\right]+t\left[C_j-C_{j-2}\right]-(j-1) L_{j-2}+2, \ \textup{ if } \ j=2 ,\\
			-\left[L_j-L_{j-2}\right]+t\left[C_j-C_{j-2}\right]-(j-1) L_{j-2}, \qquad \textup{ otherwise}.
		\end{array}\right.
	\end{aligned}
\end{align*}
It implies for $j=2$,
\begin{align*}
	(j+1) L_j=t\left[C_j-C_{j-2}\right]+2,
\end{align*}
and for $j>2$
\begin{align*}
	(j+1) L_j=-(j-2) L_{j-2}+t\left[C_j-C_{j-2}\right].
\end{align*}
Consequently, we are able to evaluate all the $C_j$ and $L_j$ for $j=0,1,\cdots, n-1$ based on these recursive formulae and therefore construct the modified singular quadrature in the Nystr\"om discretization \eqref{Nystrom}. 


\subsection{Rrcursively compressed inverse preconditioning}
RCIP was proposed in \cite{helsing2013solving_RCIPtutorial} as a kernel-independent numerical algorithm used to solve the second kind boundary integral equations. Here we only give a brief introduction of RCIP, as the details are given in~\cite{helsing2013solving_RCIPtutorial}. The main idea of RCIP is to use the interpolation matrix to construct an inverse preconditioner that transforms the information from fine grid to coarse grid, which leads to a linear system on the coarse mesh that is much easier to solve.  The method is purely algebraic and independent of singular kernels.

Consider the integral equation of the second kind on the boundary $\Gamma$ with some corner points,
\begin{equation}
	(\mathcal{I} + \mathcal{K}) \bm{\phi}(x) = \bm{g}(x)  , \qquad x\in \Gamma, 
\end{equation}
where $\mathcal{K}$ is presumably a compact operator away from the corner. This is not mandatory in the numerical practice, as in the elastic scattering case, the double layer boundary operator and its adjoint are not compact even on the smooth boundaries~\cite{hsiao2008boundary}.  According to the mesh generation illustrated by Figure \ref{meshes}, we denote the four coarse grid panels near each corner point by $\Gamma^*_k$, $\ k = 1,2,\cdots,M$, where $M$ is the number of corner points. We then split the operator into $\mathcal{K}=\mathcal{K}^* + \mathcal{K}^o$, where $\mathcal{K}^*$ is the operator when the target point $\bx$ and source point $\by$ are both at the same $\Gamma^*_k$. Thus the equation is formulated into the following linear system of equations on the coarse and fine meshes
\begin{align}
	(\mathbb{I}_{coa} + \mathbb{K}^*_{coa} + \mathbb{K}^o_{coa}) \bm{\phi}_{coa} = \mathbf{g}_{coa},\label{cor_mesh_ls}\\
	(\mathbb{I}_{fin} + \mathbb{K}^*_{fin} + \mathbb{K}^o_{fin}) \bm{\phi}_{fin} = \mathbf{g}_{fin}. \label{fin_mesh_ls}
\end{align}
Here $\mathbb{I}$ is the identity matrix, $\mathbb{K}^*$ and $\mathbb{K}^o$ are square matrices, $\{\bm{\phi}_{coa},\bm{\phi}_{fin}\}$ and $\{\mathbf{g}_{coa},\mathbf{g}_{fin}\}$ are the discretized densities and incident waves on the coarse and fine meshes respectively. 

Next, we introduce the prolongation matrix $\mathbb{P}$, which maps points on the coarse grid to points on the fine grid by Gauss-Legendre polynomial interpolation. It holds
\begin{equation}
	\mathbf{g}_{fin} \approx \mathbb{P}\mathbf{g}_{coa}.
\end{equation}
The weighted matrix $\mathbb{P}_W = \mathbb{W}_{fin}\mathbb{P}\mathbb{W}_{coa}^{-1}$ also needs to be constructed, where $\mathbb{W}_{coa}$ and $\mathbb{W}_{fin}$ are two diagonal matrices whose diagonal elements are the integral weights on the coarse and fine meshes respectively.  It also holds the property that 
\begin{equation}
	\mathbb{K}^o_{fin} \approx \mathbb{P}\mathbb{K}^o_{coa}\mathbb{P}_W^T.
\end{equation}
Let us introduce the substitution of variable 
\begin{align}\label{RCIPsubstitution}
	\bm{\bm{\phi}}_{fin} = (\mathbb{I}_{fin} + \mathbb{K}^*_{fin})^{-1} \tilde{\bm{\bm{\phi}}}_{fin} = (\mathbb{I}_{fin} + \mathbb{K}^*_{fin})^{-1} \mathbb{P}\tilde{\bm{\bm{\phi}}}_{coa}.
\end{align}
Substituting equation \eqref{RCIPsubstitution} into equation \eqref{fin_mesh_ls}, and applying the weighted prolongation matrix $\mathbb{P}_W$, we obtain
\begin{align}\label{RCIP}
	(\mathbb{I}_{coa} + \mathbb{K}_{coa}^o \mathbb{P}_W^T(\mathbb{I}_{fin} + \mathbb{K}^*_{fin})^{-1} \mathbb{P}) \tilde{\bm{\bm{\phi}}}_{coa} = \mathbf{g}_{coa},
\end{align}
where the equal sign is given up to a negligible numerical error. Denote $\mathbb{R} = \mathbb{P}_W^T(\mathbb{I}_{fin} + \mathbb{K}^*_{fin})^{-1} \mathbb{P}$ the compressed weighted inverse, which is complicated to evaluate directly but in practice can be evaluated through recursive operations~\cite{helsing2013solving_RCIPtutorial}. 

Therefore, the problem of solving linear equations \eqref{fin_mesh_ls} on the fine grid is transformed to the equation \eqref{RCIP} on the coarse grid,  which not only improves the numerical stability, but also greatly reduces the computational cost. Once equation \eqref{RCIP} is solved, we can obtain the solution $\hat{\bm{\bm{\phi}}}_{coa}$ for the density function on the coarse grid via
\begin{equation}
	\hat{\bm{\phi}}_{coa} = \mathbb{R} \tilde{\bm{\phi}}_{coa}.
\end{equation}
Compared to $\bm{\phi}_{coa}$ in \eqref{cor_mesh_ls}, $\hat{\bm{\phi}}_{coa}$ carries the information from the fine grid and is sufficient for the field evaluation away from the corner points. If the field near the corners is needed, one can also reconstruct $\bm{\phi}_{fin}$ on  $\Gamma_k^{\star}$  by a reverse recursion formula. More details can be found in  \cite{helsing2013solving_RCIPtutorial}. 

\subsection{Convergence}
It is generally difficult to analyze the convergence of Nystr\"om method in the presence of corners. The convergence analysis when discretizing the second kind boundary integral equation for Helmholtz equations in cornered domains can be found at~\cite{Kress1990ANM}. The proof can not be directly extended to the elastic equation due to the non-compactness of the boundary operator $\mathcal{D}$ and $\mathcal{D}'$ in equations \eqref{DND} and \eqref{SNN}. In case the boundary is smooth, we have shown in~\cite{Dong2020AHA} that the Nystr\"om discretization is convergent based on the Helmholtz decomposition and trigonometric interpolation. One of the key ingredients in the proof is to construct an appropriate regularizer for the boundary operator so that the equation is converted to the form of an invertible operator plus a compact one. One might follow the same idea with a detailed spectral analysis for boundary operators near the corner to show that the convergence of the proposed method. In this paper, we mainly focus on the numerical discretization of the elastic boundary integral equations and leave the convergence analysis in the future.

\section{Numerical experiments}\label{numer_exp}

\begin{figure}[htb]
	\centering
	\subcaptionbox{\centering Circle.}{
		\includegraphics[width=0.4\textwidth]{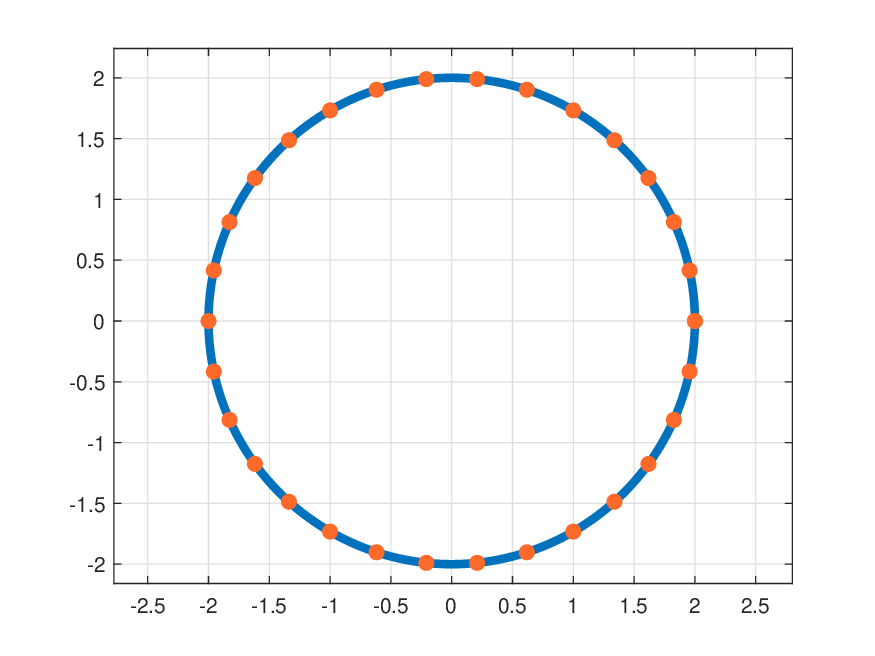}\quad
	}
	\subcaptionbox{\centering Ellipse.}{
		\includegraphics[width=0.4\textwidth]{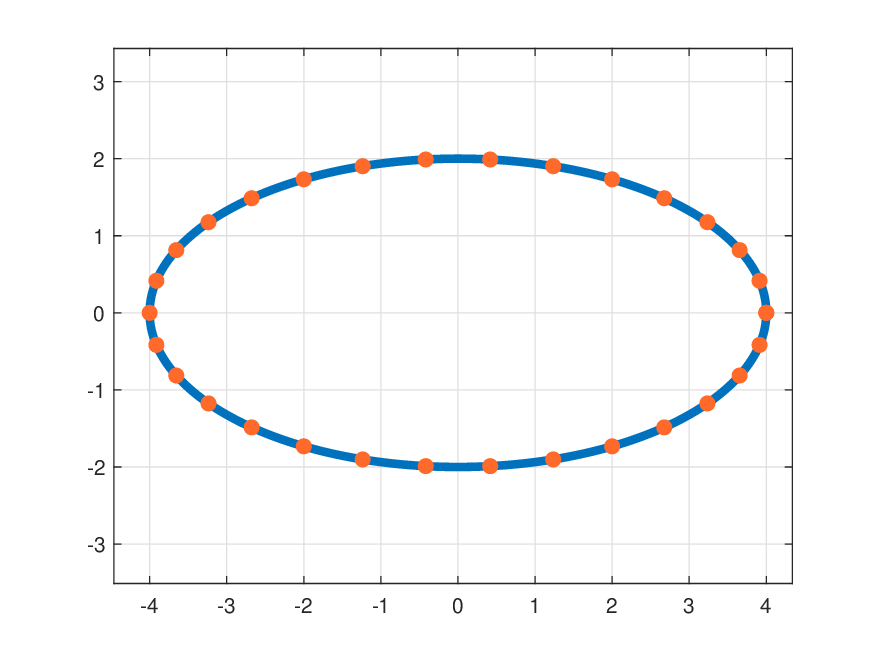}
	} 
	\\
	\subcaptionbox{\centering Droplet.}{
		\includegraphics[width=0.4\textwidth]{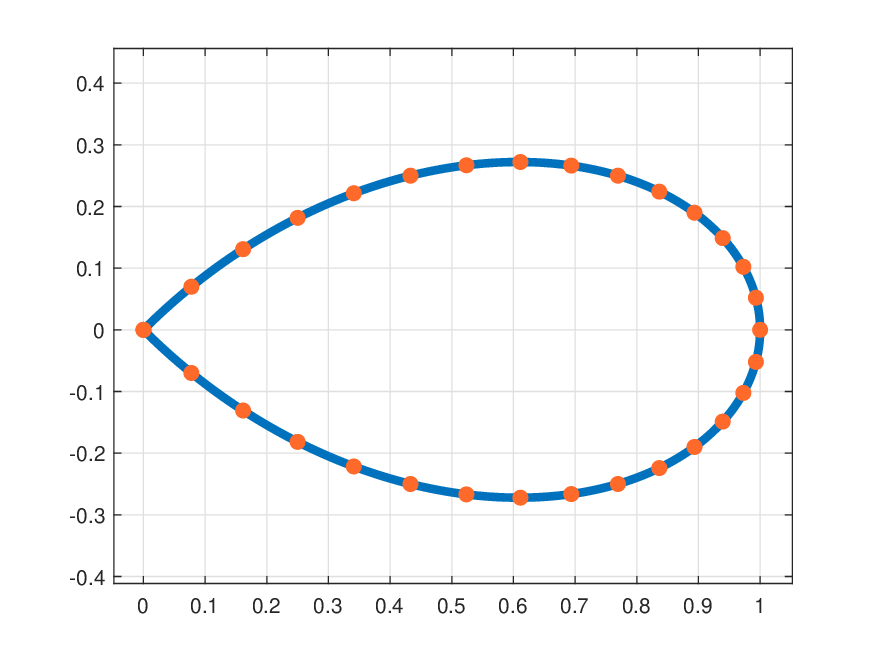}\quad
	}
	\subcaptionbox{\centering Sector.}{
		\includegraphics[width=0.4\textwidth]{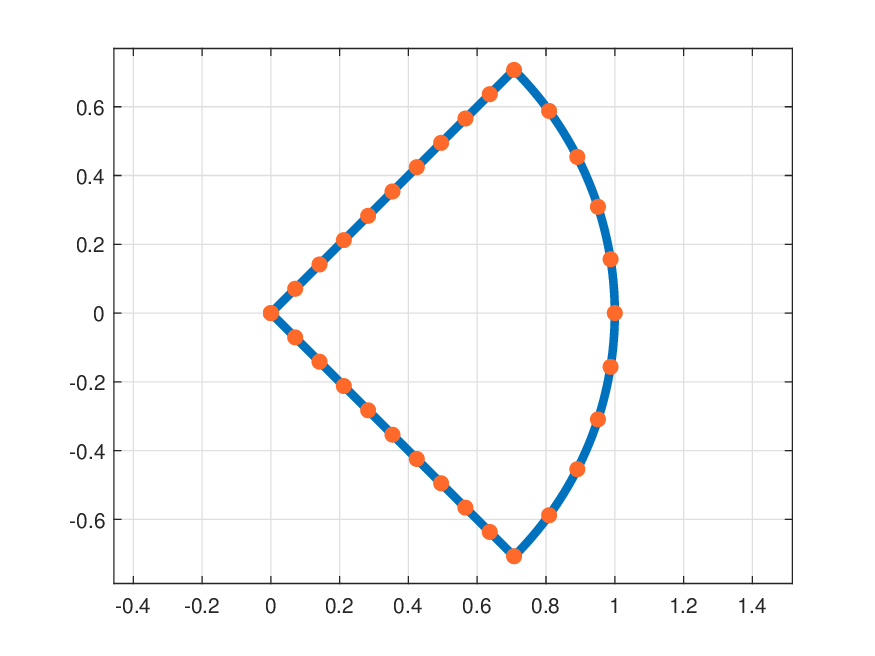}
	}
	\vspace{-2pt}
	\caption{Four different geometries in the scattering of elastic waves. } \label{geometries}
\end{figure}

In this section, we will test our algorithm by several numerical examples. The geometries used in the numerical tests are shown in Figure \ref{geometries}, namely, a circle, an ellipse, a droplet and a sector. Obviously, the droplet and the sector in Figures~\ref{geometries}(C) and \ref{geometries}(D) are non-smooth with some corner points. The prametrizations for the circle and the ellipse are given by
\begin{align}
	\bx(s) = \left\{\begin{array}{l}
		x_{1}=2\alpha\cos(s), \\
		x_{2}=2\sin(s),
	\end{array} \quad s \in[0,2\pi],\right.
\end{align}
with $\alpha = 1$ for the circle (Figure~\ref{geometries}(A))  and $\alpha=2$ for the ellipse (Figure~\ref{geometries}(B)). The droplet is given by
\begin{align} \label{droplet curve}
	\bx(s) = \left\{\begin{array}{l}
		x_{1}=\sin(\pi s)(\cos(\pi/2(s-0.5))), \\
		x_{2}=\sin(\pi s)(\sin(\pi/2(s-0.5))),
	\end{array} \quad s \in[0,1],\right.
\end{align}
 and the sector is parameterized as
\begin{align} 
	\bx(s) = (x_1(s), x_2(s)) = \left\{\begin{array}{l} \label{sector}
		\left( \beta s , -\beta k s \right), \quad s \in[0,1), \\
		\left( \cos(2\theta s - 3\theta)  ,  \sin(2\theta s - 3\theta)\right), \quad s \in[1,2), \\
		\left( -\beta s + 3 \beta, -\beta k s +3k\beta \right), \quad s \in[2,3],
	\end{array} \right.
\end{align}
where $k$ quantifies the opening of the sector with $\beta=1/\sqrt{1+k^2}$ and $\theta = \arctan(k)$.

We use $DND$ to denote solving Dirichlet problems \eqref{DND} by the double layer potential and $SNN$ to denote  solving Neumann problems \eqref{SNN} by the single layer potential, and use the subscripts ${ex}$ and ${in}$ to  represent the exterior and interior problems, respectively. We set the La\'me parameters $\lambda=1$ and $\mu=2$, the mass density $\rho=1$ and the angular frequency $\omega=3$ for all examples except example 3 where we test the performance of the solver at different wavenumbers. We choose $n=16$ for the number of Gauss-Legendre nodes placed in each panel during the discretization.  To verify the accuracy, an artificial solution for the elastic equations is constructed by the point source technique. Namely, for the interior or exterior problems with domain $\Omega$, we can put a point source outside of $\Omega$ and impose the exact boundary condition on $\Gamma$. When the boundary integral equation is solved correctly, the point source solution should be recovered in $\Omega$ by the uniqueness theorem\cite{sneddon1980vd_uniqueness_3}. All the experiments are carried out by \textit{MATLAB} on a laptop with an Intel CPU inside.

\subsection{Example 1: Elastic scattering with uniform meshes}
In this example, We test the exterior problems from the four geometries using uniform discretization. The exact solution is constructed by using a point source located at $\by_0 = (0.5,0) \in \mathbb{R}^2 \backslash \Omega $ and we verify the field at $\bx=(12.1, 5.2)$.  The accuracy of the numerical results is shown in Table \ref{table_coarse}, where  $N$ is the number of panels and  we check the errors for both the shear and compressional wave at each given $N$.

\begin{table}[h]
	\centering
	\caption{Results of elastic scattering from four geometries by uniform discretization on the boundaries. Curves A, B, C and D represent the circle, ellipse, droplet and sector, respectively. }
	\begin{tabular}{ccccccccc} 
		\hline
		\multirow{2}[4]{*}{N} & \multicolumn{4}{c}{$DND_{ex}$}     & \multicolumn{4}{c}{$SNN_{ex}$} \\
		\hline  & Curve A & Curve B & Curve C & Curve D & Curve A & Curve B & Curve C & Curve D \\
		\hline \multirow{2}[2]{*}{12} & 3.53E-17 & 5.22E-14 & 7.79E-06 & 2.96E-09 & 7.24E-17 & 9.94E-14 & 5.38E-07 & 3.45E-08 \\
		& 4.54E-17 & 4.98E-14 & 2.43E-05 & 3.37E-08 & 9.78E-17 & 2.69E-13 & 7.32E-07 & 5.27E-08 \\
	   \hline \multirow{2}[2]{*}{24} & 8.07E-17 & 1.43E-17 & 5.44E-06 & 1.41E-09 & 1.47E-16 & 2.33E-16 & 3.00E-07 & 9.52E-09 \\
		& 1.12E-16 & 4.05E-17 & 1.59E-05 & 1.19E-08 & 9.41E-17 & 5.79E-16 & 4.36E-07 & 1.41E-08 \\
		 \hline  \multirow{2}[2]{*}{36} & 2.21E-16 & 4.94E-17 & 4.39E-06 & 9.04E-10 & 3.69E-16 & 1.22E-15 & 2.14E-07 & 4.57E-09 \\
		& 3.51E-16 & 3.82E-17 & 1.24E-05 & 6.42E-09 & 2.84E-16 & 2.28E-15 & 3.23E-07 & 6.62E-09 \\
		  \hline  \multirow{2}[2]{*}{48} & 3.82E-16 & 1.34E-16 & 3.76E-06 & 6.40E-10 & 4.63E-16 & 8.80E-16 & 1.68E-07 & 2.74E-09 \\
		& 4.85E-16 & 1.30E-16 & 1.04E-05 & 4.15E-09 & 7.35E-16 & 1.46E-15 & 2.61E-07 & 3.90E-09 \\
		\hline
	\end{tabular}%
	\label{table_coarse}%
\end{table}%

As one can see from Table \ref{table_coarse}, for smooth geometries in Figures \ref{geometries}(a) and \ref{geometries}(b), an accuracy of 16 digits can be achieved by directly applying the kernel splitting method with uniform discretization. However, for boundaries with corners, the singularity at the corner causes a significant loss of accuracy. Such a loss cannot be improved by simply increasing the number of panels, as it greatly increases the number of unknowns and the condition number.  

\subsection{Example 2: Elastic scattering with non-uniform meshes}

\begin{figure}[h]
	\centering
	\begin{subfigure}{0.49\linewidth}
		\centering
		\includegraphics[width=0.9\textwidth]{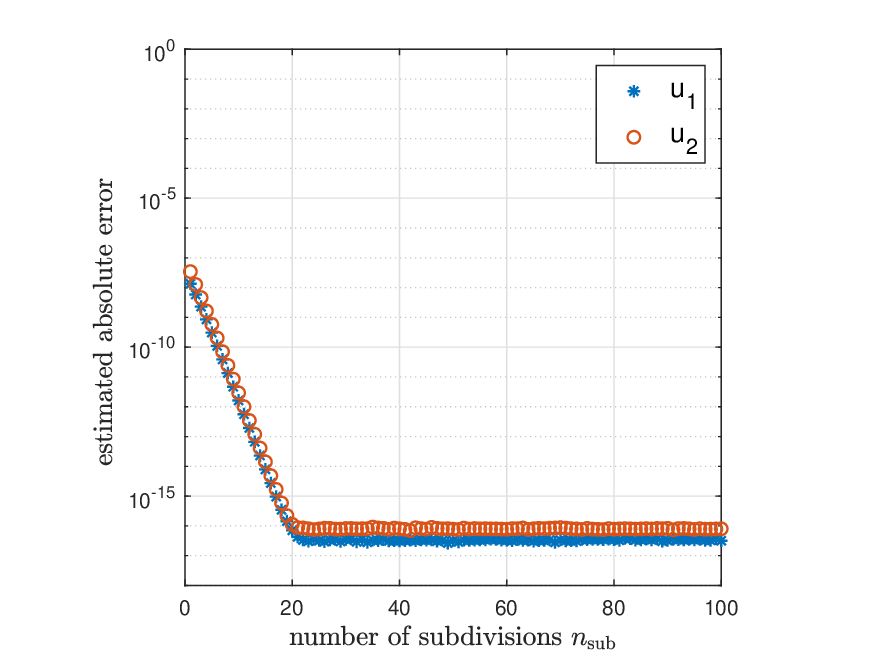}
		\caption{Neumann exterior problem}
		\label{SNN_ex}
	\end{subfigure}
	\begin{subfigure}{0.49\linewidth}
		\centering
		\includegraphics[width=0.9\textwidth]{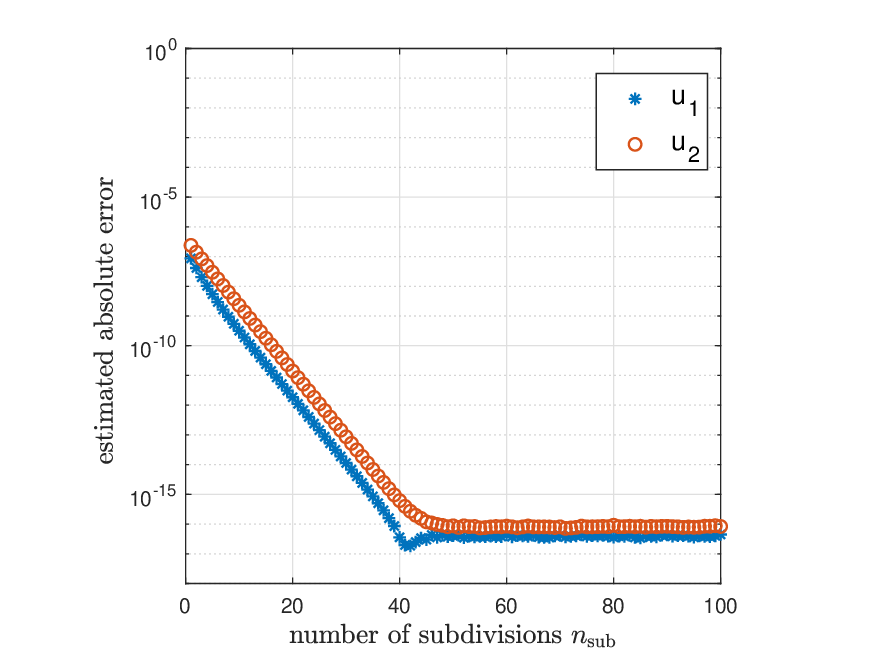}
		\caption{Dirichlet exterior problem}
		\label{DND_ex}
	\end{subfigure}
	\caption{Exterior elastic scattering problems from a droplet. A graded mesh is applied to the neighborhood of the corner and RCIP is used to solve the resulted linear system.}
	\label{nsub_convergence}
\end{figure}

In this example, we test the elastic scattering problems for the droplet by non-uniform discretization. As compared to the last example, we refine the boundary near the corner by graded mesh and compute the numerical solutions by RCIP with different refinement numbers $n_{sub}$. The results for the exterior problems are shown in Figure~\ref{nsub_convergence}.

\begin{figure}[htb]
		\centering
	\begin{subfigure}{0.49\linewidth}
		\centering
		\includegraphics[width=0.9\textwidth]{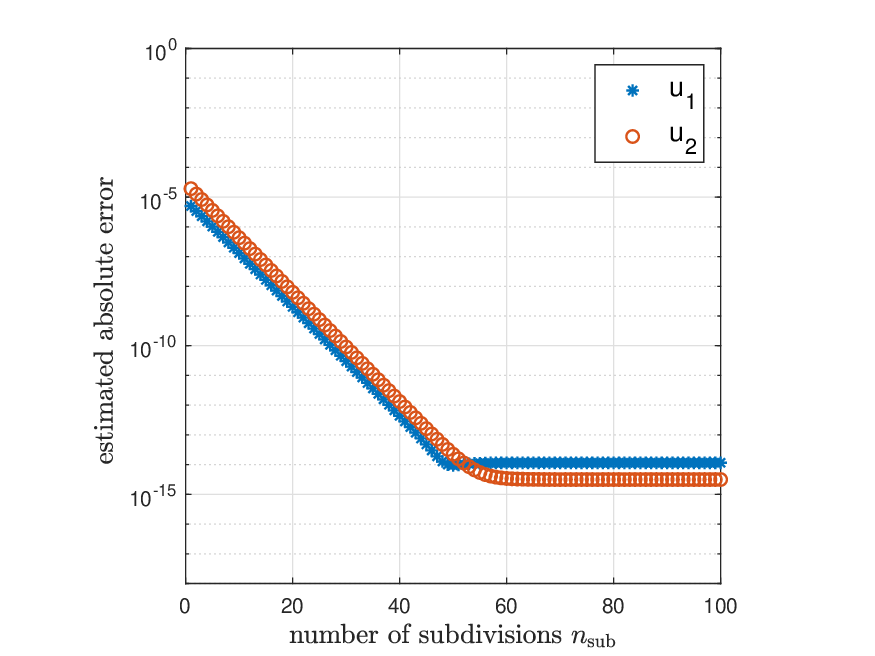}
		\caption{Neumann interior problem}
		\label{SNN_in}
	\end{subfigure}
	\begin{subfigure}{0.49\linewidth}
		\centering
		\includegraphics[width=0.9\textwidth]{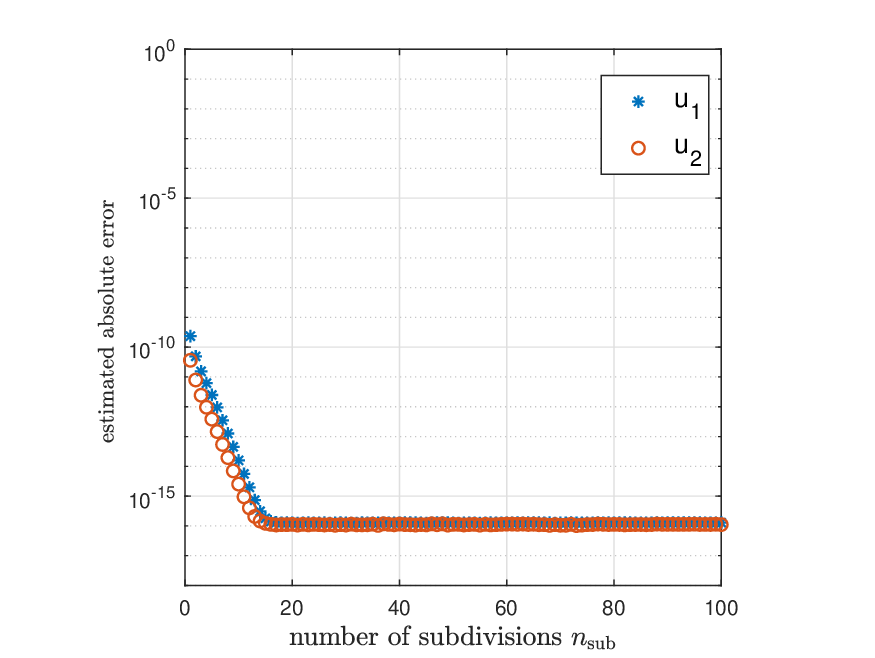}
		\caption{Dirichlet interior problem}
		\label{DND_in}
	\end{subfigure}
	\caption{Interior problems for the elastic scattering in a droplet. A graded mesh is applied to the neighborhood of the corner and RCIP is used to solve the resulted linear system.}
	\label{nsub_converge2}
\end{figure}

As one can see from Figure \ref{nsub_convergence},  with an increasing mesh refinement near the corner, the accuracy has been steadily improved and finally stabilized at the machine accuracy. One interesting fact from Figure \ref{nsub_convergence} is that the convergence rate of Neumann exterior problem is much faster than that of Dirichlet exterior problem. This is slightly counter intuitive since the solution for the Neumann problem with single layer potential is more singular than the solution for the Dirichlet problem with double layer potential~\cite{DensitySingular}. However, for the interior problems, the solution for the double layer potential does converge faster than that of the single layer potential, as shown in Figure~\ref{nsub_converge2}. In both cases, the accuracy reaches the machine precision and stays there when the number of subdivisions increases, which implies the proposed method is very stable for the resulted linear system.

\subsection{Example 3: Elastic scattering at different wavenumbers}
In this example, we test the exterior Dirichlet elastic scattering problems from a sector with different angular frequencies $\omega$, namely, from $\omega=3$ to $\omega=100$, and choose to represent the solution by the combined layer potential \eqref{combinepot} to avoid possible resonances. The corresponding numerical results are given in Table \ref{tableD+iS}. It is well known that the higher the frequency, the more oscillatory the solution becomes, which requires finer mesh to resolve the solution. In this numerical experiment, one can find that once the oscillation is well resolved, the proposed scheme reaches to the machine precision in all the frequencies that are tested.  
\begin{table}[htb]
	\centering
	\caption{Results for the elastic scattering of a sector at different wave frequency}
	\begin{tabular}{ccccccc}
		\hline
		\multirow{2}[4]{*}{N} & \multirow{2}[4]{*}{$n_{sub}$} & \multicolumn{5}{c}{Absolute Error} \\
		  &       & $\omega=3$ & $\omega=10$ & $\omega=30$ & $\omega=50$ & $\omega=100$ \\
		\hline
		\multirow{4}[2]{*}{24} & \multirow{2}[1]{*}{10} & 8.57E-14 & 1.08E-12 & 1.51E-12 & 1.42E-13 & 8.59E-09 \\
		&       & 3.60E-13 & 1.18E-12 & 6.15E-13 & 1.24E-11 & 1.02E-08 \\
		& \multirow{2}[1]{*}{20} & 2.37E-16 & 9.98E-17 & 1.11E-16 & 7.77E-12 & 8.59E-09 \\
		&       & 2.28E-16 & 4.31E-17 & 7.63E-17 & 5.02E-12 & 1.02E-08 \\
		\hline
		\multirow{4}[2]{*}{36} & \multirow{2}[1]{*}{10} & 4.62E-14 & 5.77E-13 & 8.09E-13 & 4.13E-12 & 8.20E-12 \\
		&       & 1.92E-13 & 6.30E-13 & 3.29E-13 & 6.80E-12 & 1.93E-11 \\
		& \multirow{2}[1]{*}{20} & 1.73E-16 & 1.67E-15 & 9.28E-17 & 2.89E-15 & 8.59E-12 \\
		&       & 6.13E-16 & 2.06E-16 & 6.34E-17 & 8.32E-15 & 1.84E-11 \\
		\hline
		\multirow{4}[2]{*}{72} & \multirow{2}[1]{*}{10} & 1.24E-14 & 1.97E-13 & 2.79E-13 & 1.42E-12 & 3.97E-13 \\
		&       & 6.30E-14 & 2.16E-13 & 1.13E-13 & 2.33E-12 & 1.97E-12 \\
		& \multirow{2}[1]{*}{20} & 3.56E-15 & 2.45E-15 & 8.69E-17 & 4.27E-16 & 6.82E-16 \\
		&       & 3.29E-15 & 3.07E-16 & 5.02E-17 & 1.28E-16 & 8.59E-16 \\
		\hline
	\end{tabular}
	\label{tableD+iS}
\end{table}

\begin{figure}[htb]
	\centering
	\subcaptionbox{$\bm{\phi}$ at $\omega = 3$ }{
		\includegraphics[width=0.29\textwidth]{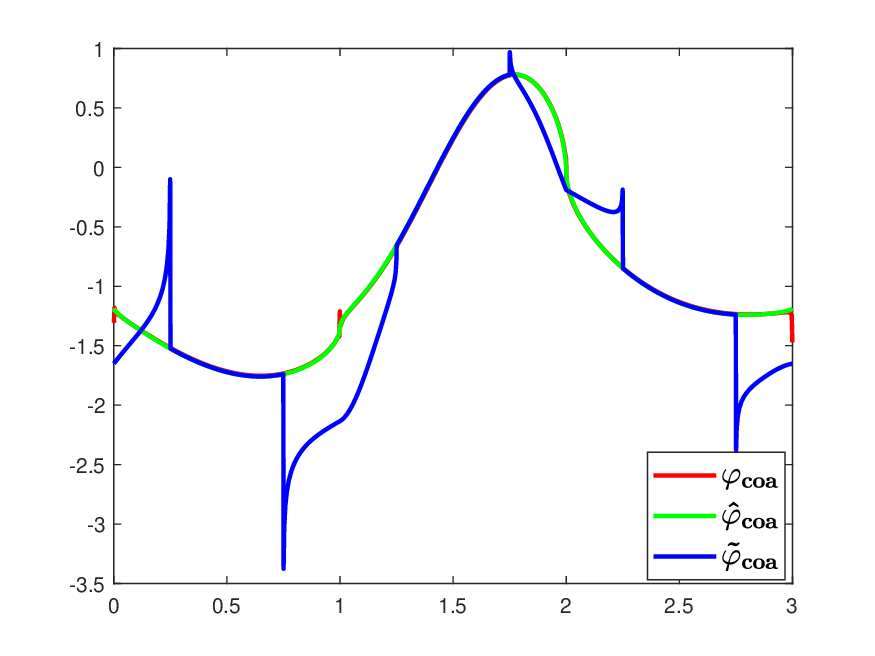}
		\quad
	}
	\subcaptionbox{$\bm{\phi}$ at $\omega = 30$}{
		\includegraphics[width=0.29\textwidth]{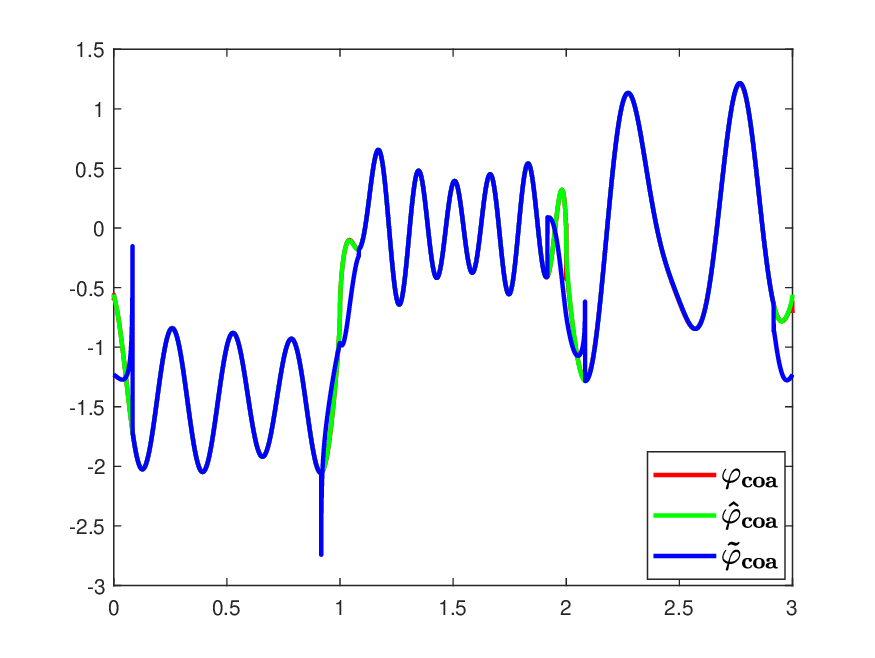}
		\quad
	} 
	\subcaptionbox{$\bm{\phi}$ at $\omega = 100$}{
		\includegraphics[width=0.29\textwidth]{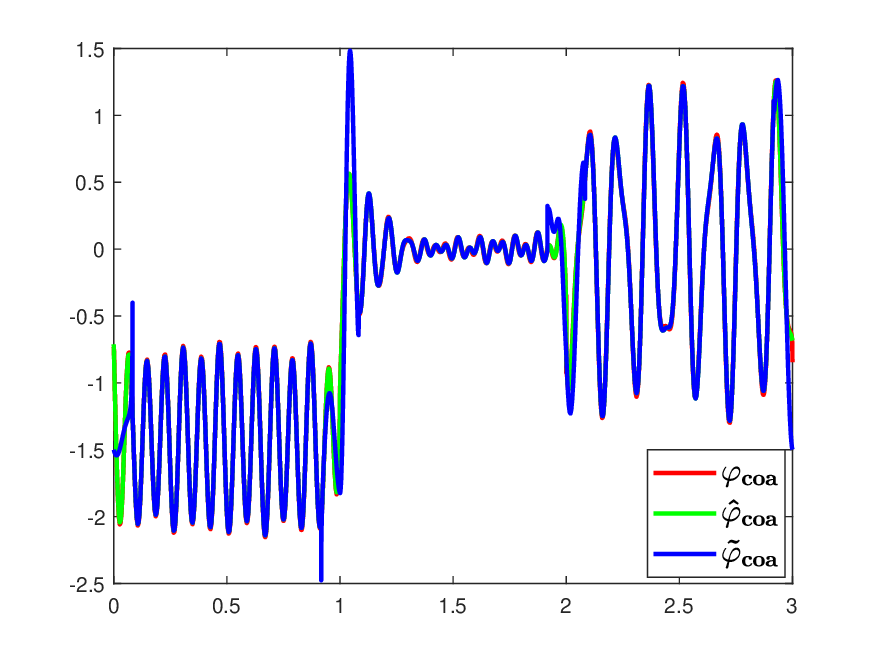}
	}
	\vspace{-2pt}
	\caption{Density functions at different wave frequencies} 
	\label{density_function} 
\end{figure}

We also consider incidence of the compressional plane wave with incident angle $d = (\cos{\pi/4},\sin{\pi/4})$. The density functions for different frequencies are shown in Figure \ref{density_function}. From the parametrization equation \eqref{sector}, the corresponding corner points are located at $s=0,1,2$ (and $3$ is the same as $0$ due to the periodicity) and one can clearly see that the solution $\bm{\phi}_{coa}$ has jumps at these points and becomes highly oscillatory when the frequency increases. 

\subsection{Example 4: The leading asymptotic behavior of the density functions}
\begin{figure}[htb]
	\centering 
	\includegraphics[width=0.45\textwidth]{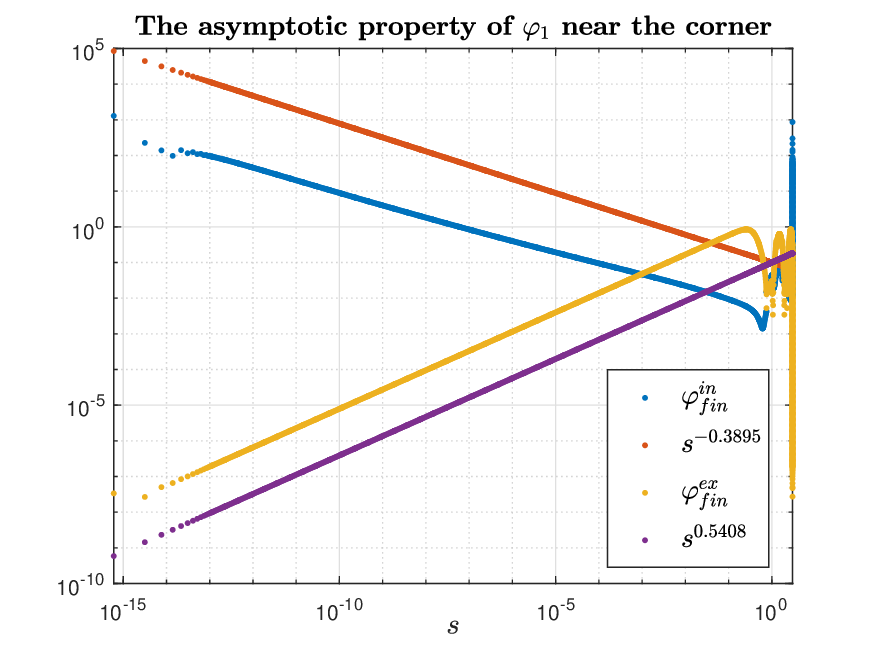}
	\includegraphics[width=0.45\textwidth]{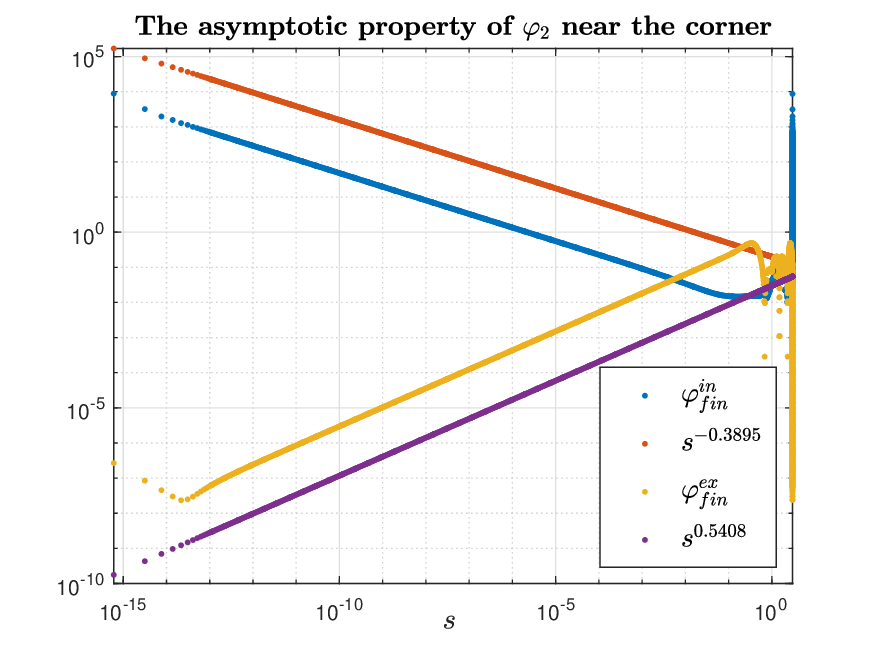}
	\vspace{-2pt}
	\caption{The asymptotic behavior of density functions in the neighborhood of the corner of a sector.} 
	\label{densityab}
\end{figure}

 As we have observed from Figure \eqref{density_function}, the density function behaves a certain singularity near the corner point \cite{DensitySingular}, which results from the power solution \eqref{powersolution} of the elastic wave equation. To have a more clear view on this singularity behavior, we evaluate the density function on the refined mesh near the corner.  The recursively compressed inverse preconditioning allows us to  recover the information on the fine grid through the backward recursion. This provides a scheme to calculate the value of the density function near the corner point and its asymptotic behavior. The detailed process of the reconstruction for $\bm{\phi}$ can be seen in~\cite[section 10]{helsing2013solving_RCIPtutorial}. Here we still consider the elastic scattering from a sector given by parametrization \eqref{sector} with an opening angle $\theta_0=\frac{\pi}{2}$. Let $r$ be the distance from the boundary to the corner point at $s=0$. The leading asymptotic behavior of the density function is
\begin{align}\label{density_asy}
	\bm{\phi}(r) \sim O(r^\alpha), \quad r \to 0,
\end{align}
where $\alpha$ is the leading asymptotic coefficient. Figure \ref{densityab} shows the asymptotic behavior of the density functions of the interior and exterior Neumann problems near the corner point. Through curve fitting, we find the asymptotic coefficients $\alpha$ in these two cases are about -0.34613 and 0.54081, respectively. 

On the other hand, by Green's representation theorem, the leading order coefficient $\alpha$ should be related to the solution of the interior and exterior problems at that corner. It was shown in \cite{NA2022} that for a corner with an opening angle $\theta$, the leading order coefficient of the elastic field $\bm{u}$ is related to the roots of some transcendental equations, which in the rigid boundary case is 
\begin{eqnarray}\label{rigidcase}
	\nu^2(1/(3-4\xi)^2)\sin^2(\theta)-\sin^2(\nu\theta)=0, \mbox{ with } \xi = \frac{\lambda}{2(\lambda+\mu)},
\end{eqnarray}
 and in the traction free case is
\begin{eqnarray}\label{tractionfree}
	\nu^2\sin^2(\theta)-\sin^2(\nu\theta)=0.
\end{eqnarray}

Since the opening angle $\theta_0$ of the corner is $\pi/2$,  the leading order coefficient should be related to $\theta=\pi/2$  and $\theta=3\pi/2$ in equations \eqref{rigidcase} and \eqref{tractionfree}. In particular, when $\theta=3\pi/2$, there is a root $\nu_1=0.61049131527757...$ for equation \eqref{rigidcase} and $\nu_2 = 0.54448373678246...$ for equation \eqref{tractionfree}. Thus we have $\alpha \approx \nu_1-1$ for the interior problem and $\alpha\approx\nu_2$ for the exterior problem, up to a small difference. This result will be useful for designing more effective numerical scheme based on the asymptotic expansion near the corner~\cite{DensitySingular}. However, we can only give a qualitative result here without detailed analysis. It is one of the future work to explore more details between the asymptotic coefficient of BIE and the opening angle of the corner.

\section{Conclusion}\label{conclusion}
In this paper, we propose an effective solver for the elastic scattering from cornered domains. The approach follows from an explicit decomposition of the kernel functions of the boundary integral operators. A high order singular quadrature scheme is proposed by the composite Gauss-Legendre nodes and modified singular integrals. The resulted linear system from cornered domain is solved based on the RCIP method. Numerical experiments show the proposed method is efficient, stable and highly accurate in various geometries and boundary conditions. Finally, we verify the leading order singular asymptotic behavior of the density functions of boundary integral operators near the corner, which is consistent with theoretical formula. Future work includes the convergence analysis of the proposed scheme and its application to the inverse elastic scattering problems.

\textbf{Acknowledgment}: The authors would like to thank Professor Shidong Jiang at Flatiron Institute for many helpful discussions.

\end{document}